\documentclass[a4paper, 11pt, twoside]{article}
\usepackage{mfpaperstuff}

\begin{document}
\renewcommand{\crefrangeconjunction}{--}
\medmuskip=3mu plus 1mu minus 1mu
\newcommand\kwo{\stackrel{\textup{\texttt{()}}}\leftarrow}
\newcommand\ob{\texttt(}
\newcommand\cb{\texttt)}
\newcommand\lfg{\scrp_{f,g}}
\newcommand\lfgp{\lfg^+}
\newcommand\lfgm{\lfg^-}
\newcommand\res{\operatorname{res}}
\newcommand\blt{\mathbf B}
\newcommand\yper[1]{\operatorname{M}^{#1}}
\newcommand\spe[1]{\operatorname{S}^{#1}}
\newcommand\mfg{\mathscr{M}^{(f,g)}}
\newcommand\ipo{i+1}
\newcommand\ipt{i+2}
\newcommand\gn{m}
\newcommand\rea[1]{\stackrel{#1}\rightharpoondown}
\newcommand\aer[1]{\stackrel{#1}\rightharpoonup}
\newcommand\hi{\hat a}
\newcommand\hj{\hat b}
\newcommand\cidt{cover-inclusive Dyck tiling}
\newcommand\cedt{cover-expansive Dyck tiling}
\newcommand\mo{{-}1}
\newcommand\one{\mathbbm{1}}
\newcommand\bsm{\begin{smallmatrix}}
\newcommand\esm{\end{smallmatrix}}
\newcommand{\rt}[1]{\rotatebox{90}{$#1$}}
\newcommand{\thmcite}[2]{\textup{\textbf{\cite[#2]{#1}}}\ }
\newcommand\zez{\mathbb{Z}/e\mathbb{Z}}
\newcommand{\hhh}{\mathcal{H}_}
\newcommand{\sect}[1]{\section{#1}}
\newcommand\cnd{(\textasteriskcentered)}
\newcommand\cne{(\textdagger)}
\newcommand\tit{Dyck tilings and the homogeneous Garnir\\ relations for graded Specht modules}
\newcommand\trinom[3]{(#1,#2,#3)!}
\newcommand\rrh[2]{#1^{[#2]}}
\newcommand\dt[2]{\scri(#1,#2)}
\renewcommand\ne[1]{\ensuremath{\mathtt{NE}(#1)}}
\newcommand\nw[1]{\ensuremath{\mathtt{NW}(#1)}}
\newcommand\se[1]{\ensuremath{\mathtt{SE}(#1)}}
\newcommand\sw[1]{\ensuremath{\mathtt{SW}(#1)}}
\newcommand\sz[1]{\ensuremath{\mathtt{S}(#1)}}
\newcommand\nz[1]{\ensuremath{\mathtt{N}(#1)}}
\newcommand\tile[1]{\ensuremath{\operatorname{tile}(#1)}}
\newcommand\st[1]{\ensuremath{\operatorname{st}(#1)}}
\newcommand\en[1]{\ensuremath{\operatorname{en}(#1)}}
\newcommand\he[1]{\ensuremath{\operatorname{ht}(#1)}}
\newcommand\de[1]{\ensuremath{\operatorname{dp}(#1)}}
\newcommand\spr[1]{{\scriptstyle\overrightarrow{\displaystyle #1}}}
\newcommand\rps[1]{{\scriptstyle\overleftarrow{\displaystyle #1}}}
\newcommand\ndt[2]{\operatorname{i}_{#1#2}}
\newcommand\et[2]{\operatorname{e}_{#1#2}}
\newcommand\fp{F}
\renewcommand\le{^{\operatorname{left}}}
\newcommand\ri{^{\operatorname{right}}}
\newcommand\leh{^{\hat{\operatorname{left}}}}
\newcommand\rih{^{\hat{\operatorname{right}}}}
\newcommand\skw[2]{\ensuremath{#1\backslash#2}}
\newcommand\xx[1]{X_{#1}}
\newcommand\xxp[1]{\xx{#1}^+}
\newcommand\xxm[1]{\xx{#1}^-}
\newcommand\rcs[3]{\partial^{#2}_{#3}#1}
\newcommand\rrs[2]{\partial^{#2}#1}
\newcommand\ccs[2]{\partial_{#2}#1}
\newcommand\ten{10}
\newcommand\eleven{11}
\newcommand\twelve{12}
\newcommand\thirteen{13}
\newcommand\fourteen{14}
\newcommand\fifteen{15}
\newcommand\sixteen{16}
\newcommand\seventeen{17}
\newcommand\eighteen{18}
\newcommand\nineteen{19}
\newcommand\twenty{20}
\newlength\sumwid\settowidth\sumwid{3cm}
\newcommand\zsum[1]{\sum_{\hbox to \sumwid{\hfil$\mathclap{\scriptstyle#1}$\hfil}}}
\Yvcentermath1\Yfrench

\title{\tit}
\msc{05E10, 20C30, 20C08}
\runninghead{Dyck tilings and the homogeneous Garnir relations}
\toptitle

\begin{abstract}
Suppose $\la$ and $\mu$ are integer partitions with $\la\supseteq\mu$. Kenyon and Wilson have introduced the notion of a \emph{\cidt} of the skew Young diagram $\skw\la\mu$, which has applications in the study of double-dimer models. We examine these tilings in more detail, giving various equivalent conditions and then proving a recurrence which we use to show that the entries of the transition matrix between two bases for a certain permutation module for the symmetric group are given by counting \cidt s. We go on to consider the inverse of this matrix, showing that its entries are determined by what we call \emph{\cedt s}. The fact that these two matrices are mutual inverses allows us to recover the main result of Kenyon and Wilson.

We then discuss the connections with recent results of Kim et al, who give a simple expression for the sum, over all $\mu$, of the number of \cidt s of $\skw\la\mu$. Our results provide a new proof of this result. Finally, we show how to use our results to obtain simpler expressions for the \emph{homogeneous Garnir relations} for the universal Specht modules introduced by Kleshchev, Mathas and Ram for the cyclotomic quiver Hecke algebras.
\end{abstract}

\section{Introduction}

The motivation for this paper is the study of the modular representation theory of the symmetric group, and more generally the representation theory of the cyclotomic Hecke algebra of type $A$. This area has recently been revolutionised by the discovery by Brundan and Kleshchev of new presentations for these algebras, which show in particular that the algebras are non-trivially $\bbz$-graded. The contribution in the present paper concerns the definition of the \emph{Specht modules}, which play a central role in the representation theory of cyclotomic Hecke algebras. These modules have been studied within the graded setting by Brundan, Kleshchev and Wang, and developed further by Kleshchev, Mathas and Ram, who have given a presentation for each Specht module with a single generator and a set of homogeneous relations. These relations include homogeneous analogues of the classical Garnir relations for the symmetric group, which allow the Specht module to be expressed as a quotient of a `row permutation module'. Although the homogeneous Garnir relations are in some sense simpler than their classical counterparts, their statement in \cite{kmr} is awkward in that the `Garnir elements' involved are given as linear combinations of expressions in the standard generators $\psi_1,\dots,\psi_{n-1}$ for the row permutation module which are not always reduced. Our main result concerning Specht modules is an expression for each Garnir relation as a linear combination of reduced expressions; this simplifies calculations with Specht modules, both theoretically and computationally.

But this result is a by-product of the main work in this paper, which is to consider tilings of skew Young diagrams by Dyck tiles. These tilings were introduced by Kenyon and Wilson, who defined in particular the notion of a \emph{cover-inclusive Dyck tiling}. They used these tilings to give a formula for the inverse of a certain matrix $M$ arising in the study of double-dimer models. We re-interpret the entries of $M$ in terms of what we call \emph{cover-expansive Dyck tilings}, and then, by proving recurrence relations for the numbers of cover-inclusive and cover-expansive Dyck tilings, we show that (sign-modified versions of) $M$ and $M^{-1}$ are in fact change-of-basis matrices for two natural bases for a certain permutation representation of the symmetric group. The fact that the two change-of-basis matrices are obviously mutually inverse provides a new proof of Kenyon and Wilson's result. Along the way we give a result showing several different equivalent conditions to the cover-inclusive condition.

In order to derive our result on Garnir relations, we then express the sum of the elements of one of our two bases in terms of the other; this involves defining a certain function $\fp$ on partitions, and proving a similar recurrence to the recurrence for \cidt s. Combining this with our results on change-of-basis matrices means that $\fp(\la)$ equals the sum over all partitions $\mu\subseteq\la$ of the number of \cidt s of $\skw\la\mu$. In fact, this had already been shown by Kim, and then by Kim, M\'esz\'aros, Panova and Wilson, verifying a conjecture of Kenyon and Wilson. As well as providing a new proof of this result, our results working directly with the function $\fp$ allow us to derive our application to Garnir relations without using Dyck tilings.

We now describe the structure of this paper. \cref{defnsec} is devoted to definitions. In \cref{cidtsec} we study \cidt s, giving equivalent conditions for cover-inclusiveness and then proving several bijective results which allow us to deduce recurrences for the number $\ndt\la\mu$ of \cidt s of $\skw\la\mu$. In \cref{cedtsec} we prove similar, though considerably simpler, recurrences for \cedt s. In \cref{ypmsec} we recall the Young permutation module $\mfg$ for the symmetric group; we define our two bases for this module, and use Dyck tilings to describe the change-of-basis matrices. We then introduce the function $\fp$ and use it to express the sum of the elements of the first basis in terms of the second, before summarising the relationship between our work and that of Kenyon, Kim, M\'esz\'aros, Panova and Wilson. Finally in \cref{klrsec} we give the motivating application of this work, introducing the Specht modules in the modern setting, and using our earlier results to give a new expression for the homogeneous Garnir relations.

\begin{acks}
The author is greatly indebted to David Speyer and Philippe Nadeau for comments on the online forum MathOverflow which inspired the author to introduce Dyck tilings into this paper.
\end{acks}

\section{Definitions}\label{defnsec}

\subsection{Partitions and Young diagrams}

As usual, a \emph{partition} is a weakly decreasing sequence $\la=(\la_1,\la_2,\dots)$ of non-negative integers with finite sum. We write this sum as $|\la|$, and say that $\la$ is a partition of $|\la|$. When writing partitions, we may group equal parts together with a superscript, and omit trailing zeroes, and we write the partition $(0,0,\dots)$ as $\varnothing$.

The \emph{Young diagram} of a partition $\la$ is the set
\[
\lset{(a,b)\in\bbn^2}{b\ls\la_a}.
\]
We may abuse notation by identifying $\la$ with its Young diagram; for example, we may write $\la\supseteq\mu$ to mean that $\la_i\gs\mu_i$ for all $i$. If $\la\supseteq\mu$, then the \emph{skew Young diagram} $\skw\la\mu$ is simply the set difference between the Young diagrams for $\la$ and $\mu$.

We draw (skew) Young diagrams as arrays of boxes in the plane, and except in the final section of this paper we use the \emph{Russian convention}, where $a$ increases from south-east to north-west, and $b$ increases from south-west to north-east. For example, the Young diagram of $\skw{(7^2,4,3,2^2)}{(2,1^2)}$ is as follows.
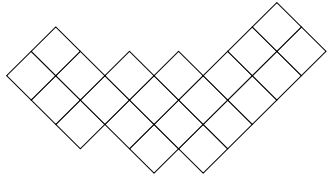
\[
\begin{tikzpicture}[scale=1,rotate=45]
\tgyoung(0cm,0cm,::;;;;;,:;;;;;;,:;;;,;;;,;;,;;)
\end{tikzpicture}
\]
The \emph{conjugate partition} to $\la$ is the partition $\la'$ obtained by reflecting the Young diagram for $\la$ left to right; thus $\la'_i=\left|\lset{j\gs1}{\la_j\gs i}\right|$ for all $i$.

We define a \emph{node} to be an element of $\bbn^2$, and a node of $\la$ to be an element of the Young diagram of $\la$. The \emph{height} of the node $(a,b)$ is $a+b$. The $j$th \emph{column} of $\bbn^2$ is the set of all nodes $(a,b)$ for which $b-a=j$.

We use compass directions to label the neighbours of a node; for example, if $\fkn$ is a node, then we write $\sw\fkn=\fkn-(0,1)$ and refer to this as the \emph{SW neighbour} of $\fkn$; we also write $\nz\fkn=\fkn+(1,1)$, and similarly for the other compass directions.

A node $\fkn$ of $\la$ is \emph{removable} if it can be removed from $\la$ to leave the Young diagram of a partition (i.e.\ if neither \nw\fkn{} nor \ne\fkn{} is a node of $\la$), while a node $\fkn$ not in $\la$ is an \emph{addable node of $\la$} if it can be added to $\la$ to leave the Young diagram of a partition.

\subsection{Tiles and tilings}

We define a \emph{tile} to be a finite non-empty set $t$ of nodes that can be ordered $\fkn_1,\dots,\fkn_r$ such that $\fkn_{i+1}\in\{\ne{\fkn_i},\se{\fkn_i}\}$ for each $i=1,\dots,r-1$. We say that $t$ \emph{starts} at its leftmost node, which we denote \st t, and \emph{ends} at its rightmost, which we denote \en t. The \emph{height} \he t of $t$ is defined to be $\max\lset{\he\fkn}{\fkn\in t}$, and we say that $t$ is a \emph{Dyck tile} if this maximum is achieved at the start and end nodes of $t$, i.e.\ $\he t=\he{\st t}=\he{\en t}$. $t$ is \emph{big} if it contains more than one node, and is a \emph{singleton} otherwise. The \emph{depth} of a node $\fkn\in t$ is
\[
\de\fkn:=\he t-\he\fkn.
\]

Now suppose $\la$ and $\mu$ are partitions with $\la\supseteq\mu$. A \emph{Dyck tiling} of $\skw\la\mu$ is a partition of $\skw\la\mu$ into Dyck tiles. Given a Dyck tiling $T$ of $\skw\la\mu$ and a node $\fkn\in\skw\la\mu$, we write \tile\fkn{} for the tile containing $\fkn$. We say that $\fkn$ is \emph{attached to} \ne\fkn{} if $\fkn$ and \ne\fkn{} lie in the same tile in $T$, and similarly for \se\fkn, \nw\fkn{} and \sw\fkn. If $t$ is a tile in $T$, the \emph{NE neighbour of $t$} is the tile starting at \ne{\en t}, if there is one, while the SW neighbour of $t$ is the tile ending at \sw{\st t}, if there is one. SE and NW neighbours of tiles are defined similarly.

Now we recast the main definition from \cite{kw}. Say that a Dyck tiling $T$ of $\skw\la\mu$ is \emph{left-cover-inclusive} if whenever $\fka$ and $\nz\fka$ are nodes of $\skw\la\mu$, \st{\tile{\nz\fka}} lies weakly to the left of \st{\tile\fka}. Similarly, $T$ is \emph{right-cover-inclusive} if whenever $\fka$ and $\nz\fka$ are nodes of $\skw\la\mu$, \en{\tile{\nz\fka}} lies weakly to the right of \en{\tile\fka}. Say that $T$ is \emph{cover-inclusive} if it is both left- and right-cover-inclusive.

We let $\dt\la\mu$ denote the set of cover-inclusive Dyck tilings of $\skw\la\mu$ if $\la\supseteq\mu$, and set $\dt\la\mu=\emptyset$ otherwise. Let $\ndt\la\mu=\left|\dt\la\mu\right|$.

\medskip
Next we make a definition which is in some sense dual to the notion of cover-inclusiveness and which appears heavily disguised in \cite{kw}. Say that a Dyck tiling is \emph{left-cover-expansive} if whenever $\fka$ and \se\fka{} are nodes of $\skw\la\mu$, \st{\tile{\se\fka}} lies weakly to the left of \st{\tile\fka}, and \emph{right-cover-expansive} if whenever $\fka$ and \sw\fka{} are nodes of $\skw\la\mu$, \en{\tile{\sw\fka}} lies weakly to the right of \en{\tile\fka}. A Dyck tiling is \emph{cover-expansive} if it is both left- and right-cover-expansive.

We write $\et\la\mu$ for the number of \cedt s of $\skw\la\mu$, setting $\et\la\mu=0$ if $\la\nsupseteq\mu$. We shall see later that $\et\la\mu\ls1$ for all $\la,\mu$.

\begin{eg}
We illustrate four Dyck tilings of $\skw{(6^2,4,3,1^2)}{(4,1^2)}$. Only the first one is cover-inclusive, and only the second is cover-expansive.
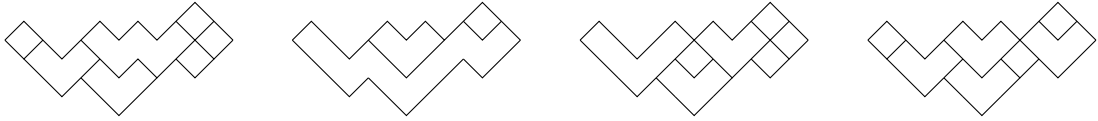
\[
\begin{tikzpicture}[scale=.25]
\draw(3,5)--++(1,-1)--++(1,1);
\draw(4,6)--++(1,-1)--++(1,1);
\draw(3,7)--++(1,-1)--++(1,1);
\draw(-6,6)--++(1,-1)--++(1,1);
\draw(-2,4)--++(1,-1)--++(1,-1)--++(1,1)--++(1,1);
\draw(-2,6)--++(1,-1)--++(1,-1)--++(1,1)--++(1,-1)--++(1,1)--++(1,1);
\draw(-5,5)--++(1,-1)--++(1,-1)--++(1,1)--++(1,1);
\draw(6,6)--++(-1,1)--++(0,0)--++(-1,1)--++(-2,-2)--++(-1,1)--++(-1,-1)--++(-1,1)--++(-2,-2)--++(-1,1)--++(0,0)--++(-1,1)
--++(-1,-1);
\end{tikzpicture}
\qquad
\begin{tikzpicture}[scale=.25]
\draw(3,5)--++(1,-1)--++(1,1);
\draw(4,6)++(1,-1)--++(1,1);
\draw(3,7)--++(1,-1)--++(1,1);
\draw(-6,6)--++(1,-1);
\draw(-2,4)--++(1,-1)--++(1,-1)--++(3,3);
\draw(-2,6)--++(1,-1)--++(1,-1)--++(2,2);
\draw(-5,5)--++(1,-1)--++(1,-1)--++(1,1);
\draw(6,6)--++(-1,1)--++(0,0)--++(-1,1)--++(-2,-2)--++(-1,1)--++(-1,-1)--++(-1,1)--++(-2,-2)--++(-1,1)--++(0,0)--++(-1,1)
--++(-1,-1);
\end{tikzpicture}
\qquad
\begin{tikzpicture}[scale=.25]
\draw(3,5)--++(1,-1)--++(1,1);
\draw(4,6)--++(1,-1)--++(1,1);
\draw(-1,3)--++(1,-1)--++(1,1);
\draw(0,4)++(1,-1)--++(1,1);
\draw(1,5)--++(1,-1)--++(1,1);
\draw(2,6)++(1,-1)--++(1,1);
\draw(3,7)--++(1,-1)--++(1,1);
\draw(-2,4)--++(1,-1);
\draw(-1,5)--++(1,-1)--++(1,1);
\draw(0,6)--++(1,-1);
\draw(-4,4)--++(1,-1)--++(1,1);
\draw(-3,5)++(1,-1)--++(1,1);
\draw(-2,6)++(1,-1)--++(1,1);
\draw(-5,5)--++(1,-1);
\draw(-6,6)--++(1,-1);
\draw(6,6)--++(-1,1)--++(0,0)--++(-1,1)--++(-2,-2)--++(-1,1)--++(-1,-1)--++(-1,1)--++(-2,-2)--++(-1,1)--++(0,0)--++(-1,1)
--++(-1,-1);
\end{tikzpicture}
\qquad
\begin{tikzpicture}[scale=.25]
\draw(3,5)--++(1,-1)--++(1,1);
\draw(4,6)++(1,-1)--++(1,1);
\draw(-1,3)--++(1,-1)--++(1,1);
\draw(0,4)++(1,-1)--++(1,1);
\draw(1,5)--++(1,-1)--++(1,1);
\draw(2,6)--++(1,-1);
\draw(3,7)--++(1,-1)--++(1,1);
\draw(-2,4)--++(1,-1);
\draw(-1,5)--++(1,-1)--++(1,1);
\draw(0,6)++(1,-1)--++(1,1);
\draw(-4,4)--++(1,-1)--++(1,1);
\draw(-3,5)++(1,-1)--++(1,1);
\draw(-2,6)--++(1,-1);
\draw(-5,5)--++(1,-1);
\draw(-6,6)--++(1,-1)--++(1,1);
\draw(6,6)--++(-1,1)--++(0,0)--++(-1,1)--++(-2,-2)--++(-1,1)--++(-1,-1)--++(-1,1)--++(-2,-2)--++(-1,1)--++(0,0)--++(-1,1)
--++(-1,-1);
\end{tikzpicture}
\]
\end{eg}

We end this section with some notation which we shall use repeatedly later. Suppose $j\in\bbz$ is fixed, and $\la$ is a partition with an addable node $\fkl$ in column $j$. We define $\xx\la$ to be the set of all integers $x$ such that $\la$ has a removable node $\fkn$ in column $j+x$, with $\he\fkn=\he\fkl-1$, and $\he\fkp<\he\fkl$ for all nodes $\fkp\in\la$ in all columns between $j$ and $j+x$. We set $\xxp\la$ to be the set of positive elements of $\xx\la$. Note that $x\in\xxp\la$ precisely when there is a Dyck tile $t\subset\la$, starting in column $j+1$ and ending in column $j+x$, which can be removed from $\la$ to leave a smaller partition; we denote this smaller partition $\rrh\la x$. We define $\rrh\la x$ for $x\in\xxm\la=\xx\la\setminus\xxp\la$ similarly.

\begin{eg}
Take $\la=(6,4^2,3,2^2)$. Then $\la$ has an addable node in column $0$, and for this node we have $\xx\la=\{-1,1,5\}$. The partitions $\rrh\la x$ are as follows.
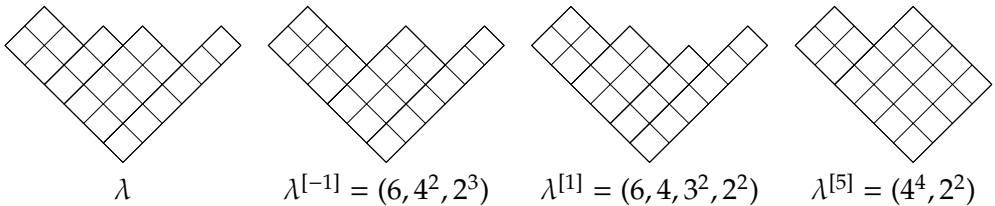
\[
\begin{array}{ccccc}
\tikz[rotate=45,scale=.8]{\tyng(0cm,0cm,6,4^2,3,2^2)}&
\tikz[rotate=45,scale=.8]{\tyng(0cm,0cm,6,4^2,2^3)}&
\tikz[rotate=45,scale=.8]{\tyng(0cm,0cm,6,4,3^2,2^2)}&
\tikz[rotate=45,scale=.8]{\tyng(0cm,0cm,4^4,2^2)}\\
\la&
\rrh\la{-1}=(6,4^2,2^3)&
\rrh\la1=(6,4,3^2,2^2)&
\rrh\la5=(4^4,2^2)
\end{array}
\]
\end{eg}

\section{Cover-inclusive Dyck tilings}\label{cidtsec}

In this section we examine \cidt s in more detail. We give some equivalent conditions to the cover-inclusive condition, and then we prove a recurrence for the number $\ndt\la\mu$ of \cidt s.

\subsection{Equivalent conditions}

\begin{thm}\label{citfae}
Suppose $\la$ and $\mu$ are partitions with $\la\supseteq\mu$, and $T$ is a Dyck tiling of $\skw\la\mu$. The following are equivalent.
\begin{enumerate}
\item\label{:ci}
$T$ is cover-inclusive.
\item\label{:de}
If $\fka$ and $\nz\fka$ are nodes of $\skw\la\mu$, then $\de{\nz\fka}\gs\de\fka$.
\item\label{:11}
If $\fka$ and $\nz\fka$ are nodes of $\skw\la\mu$, then $\tile\fka+(1,1)\subseteq\tile{\nz\fka}$.
\item\label{:nw}
If $\fka$ and $\nz\fka$ are nodes of $\skw\la\mu$ and $\fka$ is attached to $\nw\fka$, then \nw{\nz\fka} is a node of $\skw\la\mu$ and $\nz\fka$ is attached to \nw{\nz\fka}.
\item\label{:en}
If $\fka$ and $\nz\fka$ are nodes of $\skw\la\mu$ and \nz\fka{} is the end node of its tile, then $\fka$ is the end node of its tile.
\item\label{:rci}
$T$ is right-cover-inclusive.
\item\label{:ne}
If $\fka$ and $\nz\fka$ are nodes of $\skw\la\mu$ and $\fka$ is attached to $\ne\fka$, then \ne{\nz\fka} is a node of $\skw\la\mu$ and $\nz\fka$ is attached to \ne{\nz\fka}.
\item\label{:st}
If $\fka$ and $\nz\fka$ are nodes of $\skw\la\mu$ and \nz\fka{} is the start node of its tile, then $\fka$ is the start node of its tile.
\item\label{:lci}
$T$ is left-cover-inclusive.
\end{enumerate}
\end{thm}

\begin{pf}
\indent
\begin{description}
\vspace{-\topsep}
\item[{\rm(\ref{:ci}$\Rightarrow$\ref{:de})}]
Take $\fka\in\skw\la\mu$ such that $\nz\fka\in\skw\la\mu$, and let $\fkc=\st{\tile\fka}$. Since $T$ is cover-inclusive, there is a node $\fkb$ in \tile{\nz\fka} in the same column as $\fkc$, and
\[
\de{\nz\fka}\gs\he\fkb-\he{\nz\fka}\gs\he\fkc+2-(\he\fka+2)=\de\fka.
\]
\item[{\rm(\ref{:de}$\Rightarrow$\ref{:11})}]
Suppose (\ref{:11}) is false, and take $\fka\in\skw\la\mu$ such that $\nz\fka\in\skw\la\mu$ and $\tile\fka+(1,1)\nsubseteq\tile{\nz\fka}$. There is a node $\fkb\in\tile\fka$ such that $\fkb+(1,1)\notin\tile{\nz\fka}$, and we may assume that $\fkb$ is attached to $\fka$; in fact, by symmetry, we may assume $\fkb$ is either \ne\fka{} or \se\fka. If $\fkb=\ne\fka$, then \nz\fka{} is attached to neither its NE nor its SE neighbour, so is the end node of its tile, and in particular $\de{\nz\fka}=0$. On the other hand $\de\fka>\de\fkb\gs0$, contradicting (\ref{:de}). If instead $\fkb=\se\fka$, then \ne\fka{} is attached to neither its NW nor its SW neighbour, so is the start node of its tile, and has depth $0$; but $\de\fkb>\de\fka\gs0$, and again (\ref{:de}) is contradicted.
\item[{\rm(\ref{:11}$\Rightarrow$\ref{:nw})}]
This is trivial.
\item[{\rm(\ref{:nw}$\Rightarrow$\ref{:en})}]
Suppose (\ref{:en}) is false, and take $\fka\in\skw\la\mu$ as far to the left as possible such that $\nz\fka\in\skw\la\mu$ and is the end node of its tile, while $\fka$ is not the end node of its tile.

Since $\fka$ is not the end node of its tile, it is attached to either \se\fka{} or \ne\fka. If $\fka$ is attached to \se\fka, then \se\fka{} is attached to its NW neighbour, but $\nz{\se\fka}=\ne\fka{}$ is not attached to its NW neighbour, contradicting (\ref{:nw}).

So assume that $\fka$ is attached to \ne\fka. This implies in particular that $\fka$ has positive depth.

\clamn1
If $\fkc\in\tile\fka$ and $\nz\fkc\in\tile{\nz\fka}$, then $\fkc$ is not the start node of \tile\fka.
\prof
Since \tile{\nz\fka} is a Dyck tile and \nz\fka{} is its end node, we have $\he{\nz\fkc}\ls\he{\nz\fka}$. Hence $\he\fkc\ls\he\fka$, so $\de\fkc\gs\de\fka>0$.
\malc

\clamn2
For every node $\fkb\in\tile{\nz\fka}$ we have $\sz\fkb\in\tile\fka$.
\prof
If the claim is false, let $\fkb$ be the rightmost counterexample. Obviously $\fkb\neq\nz\fka$, and in particular $\fkb$ is not the end node of its tile, so $\fkb$ is attached to either \ne\fkb{} or \se\fkb. In the first case, the choice of $\fkb$ means that we have $\se\fkb\in\tile\fka$; neither $\fkb$ nor \sz\fkb{} is in \tile\fka, so \se\fkb{} is the start node of \tile\fka, contradicting Claim 1.

So we can assume $\fkb$ is attached to \se{\fkb}. The choice of $\fkb$ means that $\sz{\se\fkb}\in\tile\fka$ and is not attached to \sz\fkb. By Claim 1 \sz{\se\fkb} is not the start node of \tile\fka, so is attached to \sz{\sz\fkb}. But now \sz\fkb{} is not attached to either its NE or SE neighbour, so is the end node of its tile. \sz{\sz\fkb} is not the end node of its tile, and this contradicts the choice of $\fka$.
\malc

Now let $\fkb=\st{\tile{\nz\fka}}$. Then by Claim 2, $\fkc:=\sz\fkb\in\tile\fka$, and by Claim 1 $\fkc$ is not the start node of \tile\fka. So $\fkc$ is attached to either \nw\fkc{} or \sw\fkc. The first possibility contradicts (\ref{:nw}), since $\fkb$ is not attached to \nw\fkb, so assume that $\fkc$ is attached to \sw\fkc. But then \nw{\fkc} is attached to neither $\fkb$ nor $\fkc$, so is the end node of its tile, while \sw\fkc{} is not the end node of its tile, and this contradicts the choice of $\fka$.
\item[{\rm(\ref{:en}$\Rightarrow$\ref{:rci})}]
Suppose $T$ is not right-cover-inclusive, and take $\fka\in\skw\la\mu$ such that $\nz\fka\in\skw\la\mu$ and \tile\fka{} ends strictly to the right of \tile{\nz\fka}. Let $\fkb=\en{\tile{\nz\fka}}$; then there is a node in \tile\fka{} in the same column as $\fkb$, which we can write as $\fkb-(h,h)$ for some $h>0$. If we let $i\in\{1,\dots,h\}$ be minimal such that $\fkb-(i,i)$ is not the end node of its tile, then $\nz{\fkb-(i,i)}$ is the end node of its tile, contradicting (\ref{:en}).
\item[{\rm(\ref{:rci}$\Rightarrow$\ref{:ne})}]
Suppose $\fka,\nz\fka\in\skw\la\mu$ and $\fka$ is attached to \ne\fka. Then \en{\tile\fka} lies to the right of $\fka$, so by (\ref{:rci}) \en{\tile{\nz\fka}} does too. So \nz\fka{} is attached to either \ne{\nz\fka} or \se{\nz\fka}. But $\se{\nz\fka}=\ne\fka$ is attached to $\fka$, so \nz\fka{} is attached to \ne{\nz\fka}.
\item[{\rm(\ref{:ne}$\Rightarrow$\ref{:st})}]
This is symmetrical to the argument that \ref{:nw}$\Rightarrow$\ref{:en}.
\item[{\rm(\ref{:st}$\Rightarrow$\ref{:lci})}]
This is symmetrical to the argument that \ref{:en}$\Rightarrow$\ref{:rci}.
\item[{\rm(\ref{:lci}$\Rightarrow$\ref{:ci})}]
Since \ref{:rci}$\Rightarrow$\ref{:ne}$\Rightarrow$\ref{:st}$\Rightarrow$\ref{:lci}, right-cover-inclusive implies left-cover-inclusive. Symmetrically, left-cover-inclusive implies right-cover-inclusive, and hence cover-inclusive.
\qedhere
\end{description}
\end{pf}

We shall use these equivalent definitions of cover-inclusiveness, often without comment, in what follows. We also observe the following property of \cidt s, which we shall use without comment.

\begin{lemma}
Suppose $T$ is a \cidt{} of $\skw\la\mu$. If $\fkn\in\skw\la\mu$ is the highest node in its column in $\skw\la\mu$, then every node in \tile\fkn{} is the highest node in its column in $\skw\la\mu$.
\end{lemma}

\begin{pf}
Suppose not, and take $\fkm\in\tile\fkn$ which is not the highest node in its column. Without loss of generality we may assume $\fkm$ lies in the column immediately to the right of $\fkn$, i.e.\ $\fkm$ is either \ne\fkn{} or \se\fkn. But if $\fkm=\ne\fkn$, then the assumption $\nz\fkm\in\skw\la\mu$ means that $\nz\fkn\in\skw\la\mu$ (otherwise $\skw\la\mu$ would not be a skew Young diagram), a contradiction. So assume $\fkm=\se\fkn$. But now $\fkm$ is attached to $\nw\fkm$ and $\nz\fkm\in\skw\la\mu$, so by \cref{citfae}(\ref{:nw}) $\nz\fkn\in\skw\la\mu$, contradiction.
\end{pf}

\subsection{Recurrences}

Now we consider recurrences. We start with a simple result.

\begin{propn}\label{noadd}
Suppose $\la$ and $\mu$ are partitions and $j\in\bbz$, and that $\mu$ has an addable node $\fkm$ in column $j$, but $\la$ does not have an addable node in column $j$. Let $\mu^+=\mu\cup\{\fkm\}$. Then $\ndt\la\mu=\ndt\la{\mu^+}$.
\end{propn}

\begin{pf}
The fact that $\la$ does not have an addable node in column $j$ implies that $\la\supseteq\mu$ if and only if $\la\supseteq\mu^+$, so we may as well assume that both of these conditions hold. Given a \cidt{} of $\skw\la{\mu^+}$, we can obtain a tiling of $\skw\la\mu$ simply by adding the singleton tile $\{\fkm\}$, and it is clear that this tiling is a \cidt.

In the other direction, suppose $T$ is a \cidt{} of $\skw\la\mu$; then we claim that $\fkm$ forms a singleton tile. If we let $\fkn$ denote the highest node in column $j$ of $\skw\la\mu$, then (since $\la$ does not have an addable node in column $j$) \ne\fkn{} and \nw\fkn{} are not both nodes of $\skw\la\mu$; suppose without loss that $\ne\fkn\notin\skw\la\mu$. Then in particular $\fkn$ is not attached to \ne\fkn{} in $T$, and so (using \cref{citfae}(\ref{:ne})) $\fkm$ is not attached to \ne\fkm. $\fkm$ cannot be attached to \se\fkm{} or \sw\fkm{} (since these are not nodes of $\skw\la\mu$), and a node in a Dyck tile cannot be attached only to its NW neighbour. So $\fkm$ is not attached to any of its neighbours, i.e.\ $\{\fkm\}$ is a singleton tile as claimed. We can remove this tile, and we clearly obtain a \cidt{} of $\skw\la{\mu^+}$. So we have a bijection between $\dt\la\mu$ and $\dt\la{\mu^+}$.
\end{pf}

Now we prove a more complicated result, for which it will help us to fix some notation.

\medskip
\begin{mdframed}[innerleftmargin=3pt,innerrightmargin=3pt,innertopmargin=3pt,innerbottommargin=3pt,roundcorner=5pt,innermargin=-3pt,outermargin=-3pt]
\noindent\textbf{Notation in force for the remainder of Section \ref{cidtsec}:}

$j$ is a fixed integer, and $\la,\mu$ are partitions such that $\la$ has an addable node $\fkl$ in column $j$, while $\mu$ has an addable node $\fkm$ in column $j$. We define ${\la^+}=\la\cup\{\fkl\}$ and $\mu^+=\mu\cup\{\fkm\}$.

If $T\in\dt\la{\mu^+}$; then $\spr T$ denotes the size of the highest tile starting in column $j+1$, and $\rps T$ denotes the size of the highest tile ending in column $j-1$.
\end{mdframed}
\smallskip

Note that if $T\in\dt\la{\mu^+}$, then there must be at least one tile starting in column $j+1$, since there are more nodes in column $j+1$ of $\skw\la{\mu^+}$ than in column $j$; so $\spr T$ is well-defined, and similarly $\rps T$ is well-defined.

\medskip
To prove our main recurrence result for the numbers $\ndt\la\mu$, we give three results in which we construct bijections between sets of \cidt s. The first of these is as follows.

\begin{propn}\label{bij1}
\[
\left|\rset{T\in\dt\la{\mu^+}}{\spr T\notin\xx\la}\right|=\left|\rset{T\in\dt{\la^+}{\mu^+}}{\textup{there is a big tile in $T$ starting in column }j}\right|.
\]
\end{propn}

\begin{pf}
Note first that if $\la^+\nsupseteq\mu^+$ then $\la\nsupseteq\mu^+$, so both sides equal zero. Conversely, if $\la\nsupseteq\mu^+$, then either $\la^+\nsupseteq\mu^+$ or there are no nodes in column $j$ of $\skw{\la^+}{\mu^+}$, so again both sides are zero. So we may assume that $\la\supseteq\mu^+$. We'll construct a bijection
\[
\phi_1:\rset{T\in\dt\la{\mu^+}}{\spr T\notin\xx\la}\longrightarrow\rset{T\in\dt{\la^+}{\mu^+}}{\textup{there is a big tile in $T$ starting in column }j}.
\]
Given $T$ in the domain, let $t$ be the highest tile in $T$ starting in column $j+1$, and let $A$ be the set of nodes in column $j$ which are higher than \st t. No node in column $j+1$ can be attached to its NW neighbour, because the highest node in this column (namely \se\fkl) is not, so (by the choice of $t$) every node $\fkb$ higher than \st t in column $j+1$ is attached to \sw\fkb; hence each $\fka\in A$ is attached to \ne\fka. In a Dyck tiling a node cannot be attached only to its NE neighbour, so each $\fka\in A$ is attached to either its NW or SW neighbour. But a node in column $j-1$ cannot be attached to its NE neighbour (because the highest node in column $j-1$ is not), so every $\fka\in A$ is attached to $\nw\fka$.

Now we consider columns $x$ and $x+1$, where $x=\spr T+j=|t|+j$. Let $\fkd=\se\fkl$. Then $\fkd=\st t+(h,h)$ for some $h\gs0$, so by \cref{citfae}(\ref{:11}) \tile\fkd{} contains $t+(h,h)$. In particular, \tile\fkd\ includes the node $\fke=\en t+(h,h)$ in column $x$, and (if $\st t\neq\en t$) includes \sw\fke. Since $\fkd$ is the highest node in its column, the same is true for every node in \tile\fkd, and in particular neither \nz\fke{} nor \nw\fke{} is a node of $\la$. However, $\fke$ cannot be a removable node of $\la$, since then $\spr T=|t|$ would lie in $\xx\la$. So \ne\fke{} is a node of $\skw\la\mu$. Now \ne\fke{} is not attached to its NW neighbour (since this is $\nz\fke\notin\la$), so no node in column $x+1$ is attached to its NW neighbour, and so no node in column $x$ is attached to its SE neighbour. This implies that any node in column $x$ of positive depth must be attached to its NE neighbour, while any node in column $x$ of depth $0$ is the end node of its tile, and this tile has a NE neighbour. In particular, $t$ has a NE neighbour.

Now suppose $\fka\in A$ has depth $1$, and write $\fka=\st t+(i,i-1)$ for $i>0$. Then \tile\fka{} contains $\st t+(i,i)$, and hence contains $\en t+(i,i)$. Since $\fka$ has depth $1$, $\en t+(i,i)$ has depth $0$, and so from the last paragraph $\en t+(i,i)$ is the end node of \tile\fka, and \tile\fka{} has a NE neighbour.

Now we can construct $\phi_1(T)$, as follows.
\begin{itemize}
\item
Add the node $\fkl$ to the tiling (as a singleton tile).
\item
For each $\fka\in A$:
\begin{itemize}
\item
if $\fka$ has depth $1$ in $T$, attach \tile\fka{} to its NE neighbour;
\item
detach $\fka$ from \ne\fka{} and \nw\fka;
\item
attach $\fka$ to \se\fka;
\item
if $\fka$ has depth at least $2$ in $T$, attach \nw\fka{} to \nz\fka;
\end{itemize}
\item
Attach $t$ to is NE neighbour.
\item
Finally, attach $\fkl$ to \se\fkl.
\end{itemize}
Now let's check that $\phi_1(T)$ is a \cidt. First we check that every tile is a Dyck tile. Passing from $T$ to $\phi_1(T)$, the tiles change in the following ways:
\begin{itemize}
\item
given a tile $u\in T$ containing a node $\fka\in A$ of depth at least $2$, $\phi_1(T)$ contains the tile obtained by replacing $\fka$ with \nz\fka{} in $u$;
\item
given a tile $u$ containing a node $\fka\in A$ of depth $1$, $u$ has a NE neighbour $v$; $\phi_1(T)$ contains the tile consisting of the portion of $u$ ending at \nw\fka, and another tile obtained by combining \nz\fka, the portion of $u$ starting at \ne\fka, and $v$;
\item
let $v$ denote the NE neighbour of $t$; then $\phi_1(T)$  contains the tile obtained by joining \nw{\st t} to $t$ and $v$.
\end{itemize}
Clearly, all these new tiles are Dyck tiles. It's also clear that the last tile mentioned is a big tile starting in column $j$. So to see that $\phi_1(T)$ lies in the codomain, all that remains is to check that $\phi_1(T)$ is cover-inclusive. But this is easy using (for example) \cref{citfae}(\ref{:nw}).

As an example, we give a \cidt{} $T$ for $\la=(12^2,10,9,6^2,5,4,3^3)$, $\mu=(1)$ and $j=-1$, and its image under $\phi_1$; the nodes $\fkl$ and $\fkm$ are marked, and the nodes in $A$ are shaded and labelled with their depths.

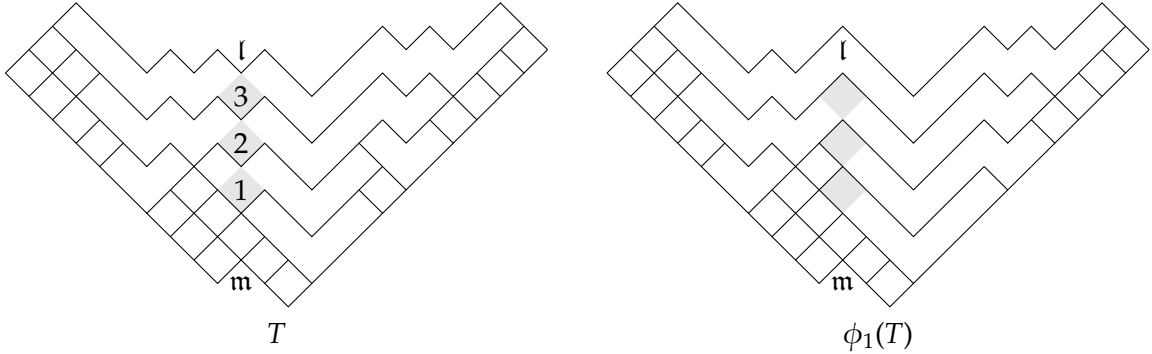
\[
\begin{array}{c@{\qquad}c}
\begin{tikzpicture}[scale=0.31]
\foreach\y in{6,8,10}\path[fill=black!10!white](-2,\y)--++(1,1)--++(1,-1)--++(-1,-1)--cycle;
\draw(0,2)--++(1,-1)--++(1,1);
\draw(4,6)--++(1,-1)--++(1,1);
\draw(5,7)--++(1,-1)--++(1,1);
\draw(6,8)++(1,-1)--++(1,1);
\draw(7,9)--++(1,-1)--++(1,1);
\draw(8,10)--++(1,-1)--++(1,1);
\draw(9,11)--++(1,-1)--++(1,1);
\draw(10,12)--++(1,-1)--++(1,1);
\draw(-1,3)--++(1,-1)--++(1,1);
\draw(4,8)--++(1,-1)++(1,1);
\draw(-3,3)--++(1,-1)--++(1,1);
\draw(-2,4)--++(1,-1)--++(1,1);
\draw(-4,4)--++(1,-1)--++(1,1);
\draw(-3,5)--++(1,-1)--++(1,1);
\draw(-5,5)--++(1,-1)--++(1,1);
\draw(-4,6)--++(1,-1)--++(1,1);
\draw(-8,8)--++(1,-1)--++(1,1);
\draw(-9,9)--++(1,-1)--++(1,1);
\draw(-8,10)--++(1,-1)++(1,1);
\draw(-10,10)--++(1,-1)--++(1,1);
\draw(-9,11)--++(1,-1)--++(1,1);
\draw(-11,11)--++(1,-1)--++(1,1);
\draw(-10,12)--++(1,-1)--++(1,1);
\draw(-1,5)--++(1,-1)--++(1,-1)--++(1,-1)--++(1,1)--++(1,1)--++(1,1);
\draw(-3,7)--++(1,-1)--++(1,-1)--++(1,1)--++(1,-1)--++(1,-1)
--++(1,1)--++(1,1)--++(1,1);
\draw(4,10)++(1,-1)--++(1,-1)--++(1,1)--++(1,1);
\draw(-9,13)--++(1,-1)--++(1,-1)--++(1,-1)--++(1,-1)--++(1,1)
--++(1,-1)--++(1,1)--++(1,-1)--++(1,1)--++(1,-1)
--++(1,-1)--++(1,1)--++(1,1)--++(1,1)--++(1,-1)
--++(1,1)--++(1,-1)--++(1,1)--++(1,1)--++(1,1);
\draw(-7,9)--++(1,-1)--++(1,-1)--++(1,1)--++(1,-1)--++(1,1)
--++(1,-1)--++(1,1)--++(1,-1)--++(1,-1)--++(1,1)
--++(1,1)--++(1,1);
\draw(-7,7)--++(1,-1)--++(1,-1)--++(1,1)--++(1,1);
\draw(12,12)--++(-1,1)--++(-1,1)--++(-2,-2)--++(-1,1)--++(-1,-1)--++(-1,1)--++(-3,-3)--++(-1,1)--++(-1,1)
--++(-1,-1)--++(-1,1)--++(-1,-1)--++(-1,1)--++(-1,-1)--++(-1,1)--++(-1,1)--++(-1,1)
--++(-3,-3);
\draw(-1,2)node{$\fkm$};
\draw(-1,12)node{$\fkl$};
\draw(-1,6)node{1};
\draw(-1,8)node{2};
\draw(-1,10)node{3};
\end{tikzpicture}
&
\begin{tikzpicture}[scale=0.31]
\foreach\y in{6,8,10}\path[fill=black!10!white](-2,\y)--++(1,1)--++(1,-1)--++(-1,-1)--cycle;
\draw(0,2)--++(1,-1)--++(1,1);
\draw(4,6)++(1,-1)--++(1,1);
\draw(5,7)--++(1,-1)--++(1,1);
\draw(6,8)++(1,-1)--++(1,1);
\draw(7,9)--++(1,-1)--++(1,1);
\draw(8,10)--++(1,-1)--++(1,1);
\draw(9,11)--++(1,-1)--++(1,1);
\draw(10,12)--++(1,-1)--++(1,1);
\draw(-1,3)--++(1,-1)--++(1,1);
\draw(-3,3)--++(1,-1)--++(1,1);
\draw(-2,4)--++(1,-1)--++(1,1);
\draw(-4,4)--++(1,-1)--++(1,1);
\draw(-3,5)--++(1,-1)--++(1,1);
\draw(-5,5)--++(1,-1)--++(1,1);
\draw(-4,6)--++(1,-1)--++(1,1);
\draw(-8,8)--++(1,-1)--++(1,1);
\draw(-9,9)--++(1,-1)--++(1,1);
\draw(-8,10)--++(1,-1)++(1,1);
\draw(-10,10)--++(1,-1)--++(1,1);
\draw(-9,11)--++(1,-1)--++(1,1);
\draw(-11,11)--++(1,-1)--++(1,1);
\draw(-10,12)--++(1,-1)--++(1,1);
\draw(-1,5)--++(1,-1)--++(1,-1)--++(1,-1)--++(1,1)--++(1,1)--++(1,1);
\draw(-3,7)--++(1,-1)--++(1,-1)++(1,1)--++(1,-1)--++(1,-1)
--++(1,1)--++(1,1)--++(1,1);
\draw(4,10)++(1,-1)--++(1,-1)--++(1,1)--++(1,1);
\draw(-9,13)--++(1,-1)--++(1,-1)--++(1,-1)--++(1,-1)--++(1,1)
--++(1,-1)--++(1,1)++(1,-1)++(1,1)--++(1,-1)
--++(1,-1)--++(1,1)--++(1,1)--++(1,1)--++(1,-1)
--++(1,1)--++(1,-1)--++(1,1)--++(1,1)--++(1,1);
\draw(-7,9)--++(1,-1)--++(1,-1)--++(1,1)--++(1,-1)--++(1,1)
--++(1,-1)++(1,1)--++(1,-1)--++(1,-1)--++(1,1)
--++(1,1)--++(1,1);
\draw(-7,7)--++(1,-1)--++(1,-1)--++(1,1)--++(1,1);
\draw(12,12)--++(-1,1)--++(-1,1)--++(-2,-2)--++(-1,1)--++(-1,-1)--++(-1,1)--++(-3,-3)--++(-1,1)--++(-1,1)
--++(-1,1)--++(-1,-1)--++(-1,-1)--++(-1,1)--++(-1,-1)--++(-1,1)--++(-1,1)--++(-1,1)
--++(-3,-3);
\draw(-2,6)--++(1,1)--++(1,-1);
\draw(-2,8)--++(1,1)--++(1,-1);
\draw(-2,10)--++(1,1)--++(1,-1);
\draw(-1,2)node{$\fkm$};
\draw(-1,12)node{$\fkl$};
\end{tikzpicture}\\
T&\phi_1(T)
\end{array}
\]

To show that $\phi_1$ is a bijection, we construct the inverse map
\[
\psi_1:\rset{T\in\dt{\la^+}{\mu^+}}{\textup{there is a big tile in $T$ starting in column }j}\longrightarrow\rset{T\in\dt\la{\mu^+}}{\spr T\notin\xx\la}.
\]
Suppose $T$ is a \cidt{} of $\skw{\la^+}{\mu^+}$ in which there is at least one big tile starting in column $j$. Since the highest node in column $j$ (namely $\fkl$) is not attached to its NE or NW neighbour, no node in column $j$ is attached to its NE or NW neighbour. So from top to bottom, column $j$ of $T$ consists of:
\begin{itemize}
\item
(possibly) some nodes attached to both their SE and SW neighbours;
\item
at least one node attached only to its SE neighbour;
\item
(possibly) some singleton tiles.
\end{itemize}

Construct $\psi_1(T)$ as follows.
\begin{itemize}
\item
For each node $\fka$ in column $j$ which is attached to \se\fka:
\begin{itemize}
\item
if $\fka$ is not attached to \sw\fka, let $\fkb$ be the first node in \tile\fka{} to the right of column $j$ which has the same height as $\fka$, and detach $\fkb$ from \sw\fkb;
\item
detach $\fka$ from \se\fka;
\item
detach $\fka$ from \sw\fka{} (if it is attached);
\item
(if $\fka\neq\fkl$) attach $\fka$ to \ne\fka{} and \nw\fka.
\end{itemize}
\item
Remove $\fkl$.
\end{itemize}
Let's check that $\psi_1(T)$ is a Dyck tiling. The tiles have changed in the following ways:
\begin{itemize}
\item
if $u\in T$ contains a node $\fka$ in column $j$ and also contains \se\fka{} and \sw\fka, then $\psi_1(T)$ contains the tile obtained by replacing $\fka$ with \sz\fka{} in $u$;
\item
if $u\in T$ is a big tile starting at a node $\fka$ in column $j$, let $\fkb$ be the first node of $u$ to the right of column $j$ which has the same height as $\fka$. Then:
\begin{itemize}
\item
$\psi_1(T)$ contains the tile consisting of the portion of $u$ starting at $\fkb$;
\item
if $u$ is the lowest big tile in $T$ starting in column $j$, then $\psi_1(T)$ also contains the tile consisting of the portion of $u$ running from $\se\fka$ to $\sw\fkb$;
\item
if $u$ is not the lowest big tile starting in column $j$, then $\psi_1(T)$ contains the tile constructed by joining together the SE neighbour of $u$, the node $\sz\fka$ and the portion of $u$ running from $\se\fka$ to $\sw\fkb$.
\end{itemize}
\end{itemize}
Again, we can see that these are all Dyck tiles, and checking that $\psi_1(T)$ is cover-inclusive is straightforward using \cref{citfae}(\ref{:nw}).

So $\psi_1(T)$ is a \cidt. It remains to check that $\spr{\psi_1(T)}\notin\xx\la$. Let $t$ be the highest tile in $\psi_1(T)$ starting in column $j+1$; then $t$ starts at \se\fka, where $\fka$ is the lowest node in column $j$ which is the start of a big tile in $T$, and ends at \sw\fkb, where $\fkb$ is the first node in \tile\fka{} to the right of column $j$ with the same height as $\fka$. In $T$, \sw\fkb{} is attached to $\fkb$, and so every node above \sw\fkb{} in the same column is attached to its NE neighbour in $T$. In particular, the highest node in this column has a NE neighbour in $\la$, and so is not removable; so $\spr{\psi_1(T)}\notin\xx\la$.

So our two maps $\phi_1,\psi_1$ really do map between the specified sets. It is easy to see from the construction that they are mutual inverses.
\end{pf}

Symmetrically, we have the following result.

\begin{propn}\label{bij1a}
\[
\left|\rset{T\in\dt\la{\mu^+}}{-\rps T\notin\xx\la}\right|=\left|\rset{T\in\dt{\la^+}{\mu^+}}{\textup{there is a big tile in $T$ ending in column }j}\right|.
\]
\end{propn}

Our next bijective result is the following.

\begin{propn}\label{bij2}
For each $x\in\xxp\la$,
\[
\left|\lset{T\in\dt\la{\mu^+}}{\spr T=x}\right|=\ndt{\rrh\la x}{\mu^+}.
\]
\end{propn}

\begin{pf}
Fix $x\in\xxp\la$. We define a bijection
\[
\phi_2:\lset{T\in\dt\la{\mu^+}}{\spr T=x}\longrightarrow\dt{\rrh\la x}{\mu^+}.
\]
Since $x\in\xxp\la$, there is a Dyck tile $\rho\subset\la$ consisting of the highest nodes in columns $j+1,\dots,j+x$ of $\la$. Now take $T\in\dt\la{\mu^+}$ with $\spr T=x$, and let $t$ be the highest tile in $T$ starting in column $j+1$. Then every northward translate of $t$ is an interval in a tile; in particular, \tile{\se\fkl} contains a translate $t+(h,h)$ for some $h\gs0$; since \se\fkl{} is the highest node in its column, the same is true for every node in \tile{\se\fkl}, so the translate $t+(h,h)$ coincides with $\rho$.

Consider the nodes in column $j$ of $\skw\la{\mu^+}$. Arguing as in the proof of \cref{bij1}, every node $\fka$ in column $j$ which is higher than $t$ is attached to \ne\fka{} and \nw\fka.

Next we consider nodes in column $j+x$. The highest node in column $j+x$ of $\skw\la{\mu^+}$ is the end node of $\rho$, and so is not attached to its NE neighbour (since this node is not in $\la$); hence no node in column $j+x$ is attached to its NE neighbour. So no node in column $j+x+1$ is attached to its SW neighbour. So from top to bottom , the nodes in column $j+x+1$ higher than \en t comprise:
\begin{itemize}
\item
(possibly) some nodes attached to their NW neighbours;
\item
(possibly) some nodes attached to neither their NW nor SW neighbours.
\end{itemize}

Now we can construct $\phi_2(T)$, as follows.
\begin{itemize}
\item
For each pair $(\fka,\fkd)$ of nodes with $\fka$ in column $j$, $\fkd$ in column $j+x+1$ and $\he\fka=\he\fkd>\he t$:
\begin{itemize}
\item
detach $\fka$ from \ne\fka;
\item
attach $\fka$ to \se\fka;
\item
if $\fkd$ is not attached to \nw\fkd, then detach $\fka$ from \nw\fka;
\item
detach $\fkd$ from \nw\fkd{} (if attached);
\item
attach $\fkd$ to \sw\fkd.
\end{itemize}
\item
Remove $\rho$.
\end{itemize}

Let's check that $\phi_2(T)$ is a Dyck tiling. The new tiles are as follows.  For every pair $(\fka,\fkd)$ as above,
\begin{itemize}
\item
if $\fkd$ is attached to \nw\fkd, then $\fka$ and $\fkd$ lie in the same tile in $T$, and $\phi_2(T)$ contains this tile but with the portion lying in columns $j+1$ to $j+x$ translated downwards;
\item
if $\fkd$ is not attached to \nw\fkd, then $\en{\tile\fka}=\nw\fkd$, and $\phi_2(T)$ contains the portion of \tile\fka{} ending at \nw\fka, and the tile constructed by joining the portion of \tile{\se\fka} between columns $j+1$ and $j+x$ to $\fka$ and \tile\fkd.
\end{itemize}
It is easy to see that these tiles are Dyck tiles, and the cover-inclusive property follows very easily from the cover-inclusive property for $T$ and the construction.

We give an example with $\la=(13^2,11,9^2,8,7,6,5^2,2^2)$, $\mu^+=(2,1)$, $s=-1$ and $x=5$. The nodes $\fkl$ and $\fkm$ are marked, and the nodes lying above $t$ in columns $s$ and $s+x+1$ are shaded.

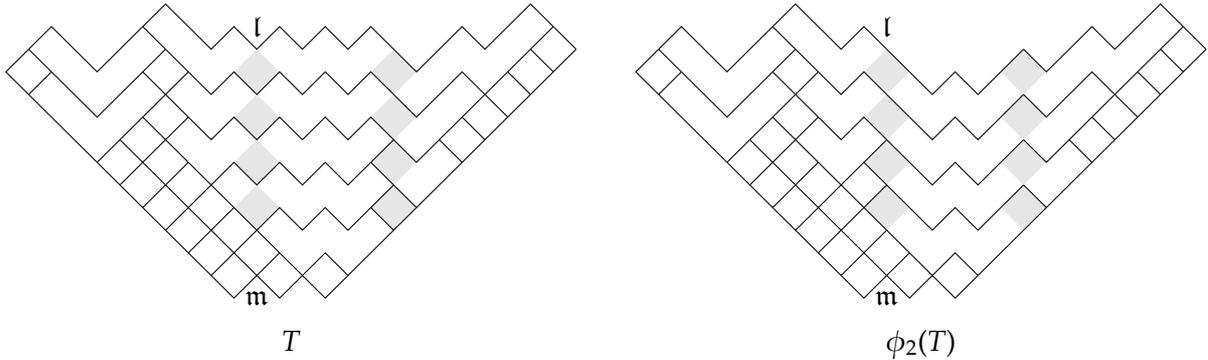
\[
\begin{array}{c@\qquad c}
\begin{tikzpicture}[scale=0.3]
\foreach\x in{-2,4}\foreach\y in{4,6,8,10}\path[fill=black!10!white](\x,\y)--++(1,1)--++(-1,1)--++(-1,-1)--cycle;
\draw(-1,3)--++(-2,-2)--++(-10,10)--++(1,1);
\foreach\x in{4,5,6,7,8,9,12}\draw(-\x,\x-2)--++(1,1);
\draw(-2,2)--++(-1,1);
\draw(0,2)--++(1,-1)--++(1,1);
\draw(-2,2)--++(1,-1)--++(1,1);
\draw(3,5)--++(1,-1)--++(1,1);
\draw(6,8)--++(1,-1)--++(1,1);
\draw(7,9)--++(1,-1)--++(1,1);
\draw(8,10)--++(1,-1)--++(1,1);
\draw(9,11)--++(1,-1)--++(1,1);
\draw(10,12)--++(1,-1)--++(1,1);
\draw(-4,4)--++(1,-1)--++(1,1);
\draw(-5,5)--++(1,-1)--++(1,1);
\draw(-6,6)--++(1,-1)--++(1,1);
\draw(-7,7)--++(1,-1)--++(1,1);
\draw(-8,8)--++(1,-1)--++(1,1);
\draw(-7,9)--++(1,-1)--++(1,1);
\draw(-9,9)++(1,-1)--++(1,1);
\draw(-8,10)++(1,-1)--++(1,1);
\draw(-7,11)--++(1,-1)--++(1,1);
\draw(-2,4)--++(1,-1)--++(1,-1)--++(1,1)--++(1,-1)--++(1,1)--++(1,1);
\draw(3,7)--++(1,-1)--++(1,-1)--++(1,1)--++(1,1);
\draw(-4,6)--++(1,-1)--++(1,-1)--++(1,1)--++(1,-1)--++(1,1)
--++(1,-1)--++(1,1)--++(1,1);
\draw(-6,10)--++(1,-1)--++(1,-1)--++(1,1)--++(1,-1)--++(1,1)
--++(1,-1)--++(1,1)--++(1,-1)--++(1,1)--++(1,-1)
--++(1,-1)--++(1,1)--++(1,1)--++(1,1);
\draw(-7,13)--++(1,-1)--++(1,-1)--++(1,-1)--++(1,1)--++(1,-1)
--++(1,1)--++(1,-1)--++(1,1)--++(1,-1)--++(1,1)
--++(1,-1)--++(1,-1)--++(1,1)--++(1,1)--++(1,-1)
--++(1,1)--++(1,1)--++(1,1);
\draw(-6,8)--++(1,-1)--++(1,-1)--++(1,1)--++(1,-1)--++(1,1)
--++(1,-1)--++(1,1)--++(1,-1)--++(1,1)--++(1,1);
\draw(-12,12)--++(1,-1)--++(1,-1)--++(1,-1)--++(1,1)--++(1,1)--++(1,1);
\draw(12,12)--++(-1,1)--++(-1,1)--++(-2,-2)--++(-1,1)--++(-2,-2)--++(-1,1)--++(-1,1)--++(-1,-1)--++(-1,1)
--++(-1,-1)--++(-1,1)--++(-1,-1)--++(-1,1)--++(-1,-1)--++(-1,1)--++(-1,1)--++(-3,-3)--++(-1,1)
--++(-1,1)--++(-1,-1);
\draw(-2,1)node{$\fkm$};
\draw(-2,13)node{$\fkl$};
\end{tikzpicture}
&
\begin{tikzpicture}[scale=0.3]
\draw(-2,1)node{$\fkm$};
\draw(-2,13)node{$\fkl$};
\foreach\x in{-2,4}\foreach\y in{4,6,8,10}\path[fill=black!10!white](\x,\y)--++(1,1)--++(-1,1)--++(-1,-1)--cycle;
\draw(0,2)--++(1,-1)--++(1,1);
\draw(-1,3)--++(-2,-2)--++(-10,10)--++(1,1);
\foreach\x in{4,5,6,7,8,9,12}\draw(-\x,\x-2)--++(1,1);
\draw(-2,2)--++(-1,1);
\draw(-2,2)--++(1,-1)--++(1,1);
\draw(3,5)++(1,-1)--++(1,1);
\draw(6,8)--++(1,-1)--++(1,1);
\draw(7,9)--++(1,-1)--++(1,1);
\draw(8,10)--++(1,-1)--++(1,1);
\draw(9,11)--++(1,-1)--++(1,1);
\draw(10,12)--++(1,-1)--++(1,1);
\draw(-4,4)--++(1,-1)--++(1,1);
\draw(-5,5)--++(1,-1)--++(1,1);
\draw(-6,6)--++(1,-1)--++(1,1);
\draw(-7,7)--++(1,-1)--++(1,1);
\draw(-8,8)--++(1,-1)--++(1,1);
\draw(-7,9)--++(1,-1)--++(1,1);
\draw(-9,9)++(1,-1)--++(1,1);
\draw(-8,10)++(1,-1)--++(1,1);
\draw(-7,11)--++(1,-1)--++(1,1);
\foreach\y in{5,7}\draw(-3,\y)--++(1,1)--++(1,-1);
\foreach\y in{10,12}\draw(-2,\y)--++(1,-1);
\draw(3,9)--++(1,1);
\draw(3,11)--++(1,1);
\draw(-2,4)--++(1,-1)--++(1,-1)--++(1,1)--++(1,-1)--++(1,1)--++(1,1);
\draw(3,7)++(1,-1)--++(1,-1)--++(1,1)--++(1,1);
\draw(-4,6)--++(1,-1)--++(1,-1)++(1,1)--++(1,-1)--++(1,1)
--++(1,-1)--++(1,1)--++(1,1);
\draw(-6,10)--++(1,-1)--++(1,-1)--++(1,1)--++(1,-1)++(1,1)
--++(1,-1)--++(1,1)--++(1,-1)--++(1,1)++(1,-1)
--++(1,-1)--++(1,1)--++(1,1)--++(1,1);
\draw(-7,13)--++(1,-1)--++(1,-1)--++(1,-1)--++(1,1)--++(1,-1)
++(1,1)--++(1,-1)--++(1,1)--++(1,-1)--++(1,1)
++(1,-1)--++(1,-1)--++(1,1)--++(1,1)--++(1,-1)
--++(1,1)--++(1,1)--++(1,1);
\draw(-6,8)--++(1,-1)--++(1,-1)--++(1,1)--++(1,-1)++(1,1)
--++(1,-1)--++(1,1)--++(1,-1)--++(1,1)--++(1,1);
\draw(-12,12)--++(1,-1)--++(1,-1)--++(1,-1)--++(1,1)--++(1,1)--++(1,1);
\draw(12,12)--++(-1,1)--++(-1,1)--++(-2,-2)--++(-1,1)--++(-2,-2)--++(-1,1)++(-6,0)
--++(-1,1)--++(-1,-1)--++(-1,1)--++(-1,1)--++(-3,-3)--++(-1,1)
--++(-1,1)--++(-1,-1);
\end{tikzpicture}\\
T&\phi_2(T)
\end{array}
\]
Now we construct the inverse map
\[
\psi_2:\dt{\rrh\la x}{\mu^+}\longrightarrow\lset{T\in\dt\la{\mu^+}}{\spr T=x}.
\]
Suppose $T$ is a \cidt{} of $\skw{\rrh\la x}{\mu^+}$. No node in column $j-1$ or column $j$ of $\rrh\la x$ can be attached to its NE neighbour (since the NE neighbours of the highest nodes in these columns are not nodes of $\rrh\la x$). Hence from top to bottom the nodes in column $j$ comprise:
\begin{itemize}
\item
(possibly) some nodes attached to their NW and SE neighbours;
\item
(possibly) some nodes attached only to their SE neighbours;
\item
(possibly) some singleton nodes.
\end{itemize}
Note that if $\fka$ is a node in column $j$ which is not attached to \nw\fka, then \nw\fka{} is the end node of its tile; in particular, \tile\fka{} has a NW neighbour.

\clam
Suppose $r\gs1$, and the node $\fka:=\fkl-(r,r)$ is attached to \se\fka. Then \tile\fka{} reaches column $j+x+1$ and includes all the nodes in $\rho-(r,r)$.
\prof
We proceed by induction on $r$. First suppose $r=1$. Since $x\in\xx\la$, every node in columns $j+1,\dots,j+x$ of $\la$ has height less than \he\fkl. Hence every node in columns $j+1,\dots,j+x$ of $\rrh\la x$ has height less than $\he\fkl-2=\he\fka$; since \tile\fka{} is a Dyck tile, it must reach \he\fka{} at some point to the right of \se\fka, and so must reach column $x+j+1$. Furthermore, $\fka$ is the highest node in its column, so every node in \tile\fka{} is the highest in its column, and in particular the portion of \tile\fka{} between columns $j+1$ and $j+x$ must consist of the highest nodes in these columns, i.e.\ the nodes in $\rho-(1,1)$.

Now suppose $r>1$. By induction \tile\fka{} cannot include any of the nodes in $\rho-(v,v)$ for any $v<r$; all the remaining nodes in columns $j+1,\dots,j+r$ have height less than \he\fka, and so (since \tile\fka{} must reach \he\fka{} at some point to the right of \se\fka{}) \tile\fka{} must reach column $j+x+1$. The translate $\tile\fka+(1,1)$ is contained in \tile{\nz\fka}, and by induction this includes the nodes in $\rho-(r-1,r-1)$; so \tile\fka{} includes the nodes in $\rho-(r,r)$.
\malc

Now we consider the nodes in column $j+x+1$. If $\fkd$ is a node in column $j+x+1$, then by similar arguments to those used above, $\fkd$ cannot be attached to \nw\fkd, and if $\fkd$ is attached to \sw\fkd, then \tile\fkd{} includes a node in column $j$ and also includes all the nodes in $\rho-(r,r)$, where $r=\frac12(\he\fkl-\he\fkd)$.

So we find that, if $\fka$ and $\fkd$ are nodes in columns $j$ and $j+x+1$ respectively with the same height, then $\fka$ is attached to \se\fka{} if and only if $\fkd$ is attached to \sw\fkd, and that in this case $\tile\fka=\tile\fkd$. Say that such a pair $(\fka,\fkd)$ is a \emph{connected} pair. Now we can construct $\psi_2(T)$ as follows.
\begin{itemize}
\item
Add $\rho$ as a tile.
\item
For each connected pair $(\fka,\fkd)$:
\begin{itemize}
\item
detach $\fkd$ from \sw\fkd;
\item
if $\fka$ is attached to \nw\fka, then attach $\fkd$ to \nw\fkd;
\item
detach $\fka$ from \se\fka;
\item
attach $\fka$ to \ne\fka;
\item
attach $\fka$ to \nw\fka{} (if it is not already attached).
\end{itemize}
\end{itemize}
To check that the resulting tiling is a Dyck tiling, we examine the new tiles:
\begin{itemize}
\item
for each connected pair $(\fka,\fkd)$ with $\fka=\fkl-(r,r)$ and $\fka$ attached to \nw\fka, $\psi_2(T)$ contains the tile obtained from \tile\fka{} by replacing the segment $\rho-(r,r)$ with $\rho-(r-1,r-1)$;
\item
for each connected pair $(\fka,\fkd)$ with $\fka=\fkl-(r,r)$ and $\fka$ not attached to \nw\fka, $\psi_2(T)$ contains the tile comprising $\fka$, the NW neighbour of \tile\fka{} and $\rho-(r-1,r-1)$, and another tile consisting of the portion of \tile\fka{} starting at $\fkd$;
\item
$\psi_2(T)$ contains the tile $\rho-(r,r)$, where $r$ is the number of connected pairs.
\end{itemize}
Once more, we can see that these are all Dyck tiles, and the cover-inclusive property is easy to check from the construction. Furthermore, the last tile mentioned above is the highest tile starting in column $j+1$, and has size $x$.

So $\psi_2$ really does map $\dt{\rrh\la x}{\mu^+}$ to $\lset{T\in\dt\la{\mu^+}}{\spr T=x}$. And it is easy to see that $\phi_2$ and $\psi_2$ are mutual inverses.
\end{pf}

Symmetrically, we have the following.

\begin{propn}\label{bij2a}
For each $x\in\xxm\la$,
\[
\left|\lset{T\in\dt\la{\mu^+}}{\rps T=-x}\right|=\ndt{\rrh\la x}{\mu^+}.
\]
\end{propn}

Our third bijection involves the partition ${\la^+}$.

\begin{propn}\label{bij3}
\[
\left|\rset{T\in\dt\la\mu}{\parbox{85pt}{\rm$\fkm$ does not lie in a singleton tile of $T$}}\right|=\left|\rset{T\in\dt{\la^+}{\mu^+}}{\parbox{129pt}{\rm there is no big tile in $T$ starting or ending in column $j$}}\right|.
\]
\end{propn}

\begin{pf}
First observe that $\la\supseteq\mu$ if and only if $\la^+\supseteq\mu^+$, so we may as well assume that both these conditions hold. Next note that if $\fkl=\fkm$, then $\skw\la\mu=\skw{\la^+}{\mu^+}$, and this skew Young diagram contains no nodes in column $j$; in particular, it does not include $\fkm$. So in this case, \cref{bij3} amounts to the trivial statement $\ndt\la\mu=\ndt{\la^+}{\mu^+}$.

So we assume that $\la\supset\mu$ and $\fkl\neq\fkm$; so $\fkm$ is a node of $\skw\la\mu$. Now we want to construct a bijection
\[
\phi_3:\lset{T\in\dt\la\mu}{\fkm\textup{ lies in a big tile in $T$ }}\longrightarrow\rset{T\in\dt{\la^+}{\mu^+}}{\parbox{129pt}{there is no big tile in $T$ starting or ending in column $j$}}.
\]

Suppose we have $T\in\dt\la\mu$ with $\fkm$ lying in a big tile. Then $\fkm$ is attached to both \ne\fkm{} and \nw\fkm, and hence every node $\fka$ in column $j$ is attached to both \ne\fka{} and \nw\fka. In particular, every node in column $j$ has depth at least $1$. Now we construct $\phi_3(T)$ from $T$ as follows.
\begin{itemize}
\item
Add $\fkl$ as a singleton tile.
\item
For each node $\fka$ in column $j$:
\begin{itemize}
\item
if $\fka$ has depth at least $2$, attach \nz\fka{} to both \ne\fka{} and \nw\fka;
\item
(if $\fka\neq\fkl$) detach $\fka$ from \ne\fka{} and \nw\fka.
\end{itemize}
\item
Remove $\fkm$.
\end{itemize}
The new tiles in $\phi_3(T)$ are as follows:
\begin{itemize}
\item
for $u\in T$ containing a node $\fka$ in column $j$ of depth at least $2$, $\phi_3(T)$ contains the tile obtained from $u$ by replacing $\fka$ with \nz\fka{};
\item
for $u\in T$ containing a node $\fka$ in column $j$ of depth $1$, $\phi_3(T)$ contains the portion of $u$ ending in column $j-1$, the portion of $u$ starting in column $j+1$, and the singleton tile $\{\nz\fka\}$.
\end{itemize}

Clearly all these tiles are Dyck tiles. The cover-inclusive property follows easily from the corresponding property of $T$, and in particular the fact (\cref{citfae}(\ref{:de})) that depth weakly decreases down columns in a cover-inclusive tiling. Moreover, $\phi_3(T)$ has no big tile starting or ending in column $j$, so lies in the codomain.

To show that $\phi_3$ is a bijection, we construct its inverse
\[
\psi_3:\rset{T\in\dt{\la^+}{\mu^+}}{\parbox{129pt}{there is no big tile in $T$ starting or ending in column $j$}}\longrightarrow\lset{T\in\dt\la\mu}{\fkm\textup{ lies in a big tile in $T$}}.
\]
Suppose $T\in\dt{\la^+}{\mu^+}$, and that there is no big tile in $T$ starting or ending in column $j$. The highest node in column $j$ of $\skw{\la^+}{\mu^+}$, namely $\fkl$, is not attached to \ne\fkl{} or \nw\fkl{} (since these are not nodes of ${\la^+}$) and so no node in column $j$ is attached to its NE or NW neighbour. So every node in column $j$ is either a singleton or attached to both its SE and SW neighbours. Now we can construct $\psi_3(T)$ from $T$ as follows.
\begin{itemize}
\item
Add $\fkm$ as singleton tile.
\item
For each node $\fka$ in column $j$:
\begin{itemize}
\item
detach $\fka$ from \se\fka{} and \sw\fka{} (if it is attached);
\item
(if $\fka\neq\fkl$) attach $\fka$ to \ne\fka{} and \nw\fka.
\end{itemize}
\item
Remove $\fkl$.
\end{itemize}
Again, it is easy to check that $\psi_3(T)$ is a \cidt, and clearly $\fkm$ lies in a big tile in $\psi_3(T)$. Furthermore, it is easy to see that $\phi_3$ and $\psi_3$ are mutually inverse.
\end{pf}

Now we combine \cref{bij1,bij1a,bij2,bij2a,bij3} to prove our main recurrence. We retain the assumptions and notation from above.

Our main result is as follows. Recall that $\ndt\la\mu$ denotes the total number of \cidt s of $\skw\la\mu$.

\begin{propn}\label{recur}
With notation as above,
\[
\ndt\la\mu+\ndt\la{\mu^+}=\ndt{\la^+}{\mu^+}+\sum_{x\in\xx\la}\ndt{\rrh\la x}{\mu^+}.
\]
\end{propn}

\begin{pf}
By \cref{bij2}, we have
\[
\left|\lset{T\in\dt\la{\mu^+}}{\spr T\in\xx\la}\right|=\sum_{x\in\xxp\la}\ndt{\rrh\la x}{\mu^+}.
\]
So
\begin{align*}
\ndt\la{\mu^+}&=\left|\lset{T\in\dt\la{\mu^+}}{\spr T\in\xx\la}\right|\ +\ \left|\lset{T\in\dt\la{\mu^+}}{\spr T\notin\xx\la}\right|\\
&=\sum_{x\in \xxp\la}\ndt{\rrh\la x}{\mu^+}\ +\ \left|\rset{T\in\dt{\la^+}{\mu^+}}{\textup{there is a big tile in $T$ starting in column }j}\right|
\end{align*}
by \cref{bij1}.

Symmetrically, we have
\[
\ndt\la{\mu^+}=\sum_{x\in\xxm\la}\ndt{\rrh\la x}{\mu^+}\ +\ \left|\rset{T\in\dt{\la^+}{\mu^+}}{\textup{there is a big tile in $T$ ending in column }j}\right|.
\]

Now consider $\ndt\la\mu$. Obviously we have
\[
\left|\lset{T\in\dt\la\mu}{\fkm\textup{ lies in a singleton tile in }T}\right|=\ndt\la{\mu^+},
\]
and by \cref{bij3}
\[
\left|\rset{T\in\dt\la\mu}{\parbox{85pt}{\rm$\fkm$ does not lie in a singleton tile of $T$}}\right|=\left|\rset{T\in\dt{\la^+}{\mu^+}}{\parbox{129pt}{\rm there is no big tile in $T$ starting or ending in column $j$}}\right|.
\]
Hence
\begin{align*}
\ndt\la\mu+\ndt\la{\mu^+}&=\left|\lset{T\in\dt\la\mu}{\fkm\textup{ lies in a singleton tile in }T}\right|+\left|\rset{T\in\dt\la\mu}{\parbox{85pt}{\rm$\fkm$ does not lie in a singleton tile of $T$}}\right|+\ndt\la{\mu^+}\\
&=\left|\rset{T\in\dt{\la^+}{\mu^+}}{\parbox{129pt}{there is no big tile in $T$ starting or ending in column $j$}}\right|+2\ndt\la{\mu^+}\\
&=\sum_{x\in\xx\la}\ndt{\rrh\la x}{\mu^+}+\left|\rset{T\in\dt{\la^+}{\mu^+}}{\textup{there is a big tile in $T$ starting in column }j}\right|\\
&\phantom{=\sum_{x\in\xx\la}\ndt{\rrh\la x}{\mu^+}\ }+\left|\rset{T\in\dt{\la^+}{\mu^+}}{\textup{there is a big tile in $T$ ending in column }j}\right|\\
&\phantom{=\sum_{x\in\xx\la}\ndt{\rrh\la x}{\mu^+}\ }+\left|\rset{T\in\dt{\la^+}{\mu^+}}{\parbox{129pt}{there is no big tile in $T$ starting or ending in column $j$}}\right|.
\end{align*}
Since a \cidt{} cannot contain big tiles starting and ending in the same column, the sum of the last three terms is $\ndt{\la^+}{\mu^+}$, and we are done.
\end{pf}

\section{Cover-expansive Dyck tilings}\label{cedtsec}

In this section we consider \cedt s, proving similar (though considerably simpler) recurrences to those in \cref{cidtsec}.

\subsection{Basic properties}

We begin by studying left- and right-\cedt s. In contrast to \cidt s, it is not the case that the left-cover-expansive and right-cover-expansive conditions are equivalent. For example, the unique Dyck tiling of $\skw{(2)}\varnothing$ is left- but not right-cover-expansive.

We begin with equivalent conditions to the left- and right-cover-expansive conditions.

\begin{propn}\label{cetfae}
Suppose $\la$ and $\mu$ are partitions with $\la\supseteq\mu$, and $T$ is a Dyck tiling of $\skw\la\mu$.
\begin{enumerate}
\item
$T$ is left-cover-expansive if and only if for every tile $t$ in $T$, we have $\nw{\st t}\notin\skw\la\mu$.
\item
$T$ is right-cover-expansive if and only if for every tile $t$ in $T$, we have $\ne{\en t}\notin\skw\la\mu$.
\end{enumerate}
\end{propn}

\begin{pf}
We prove (1); the proof of (2) is similar.

Suppose $T$ is left-cover-expansive. Given a tile $t$, let $\fka=\nw{\st t}$. If $\fka\in\skw\la\mu$, then by the left-cover-expansive property $\st{\tile\fka}$ lies weakly to the right of \st t; but $\fka$ lies strictly to the left of \st t, a contradiction.

Conversely, suppose the given property holds, and $\fka,\se\fka$ are nodes of $\skw\la\mu$. Let $\fkb=\st{\tile{\se\fka}}$, and suppose $\fkb$ lies in column $i$. Then $\nw\fkb\notin\skw\la\mu$, and hence there are no nodes in column $i-1$ of $\skw\la\mu$ higher than $\fkb$ (otherwise $\skw\la\mu$ would not be a skew Young diagram). In particular, there are no nodes of \tile\fka{} in column $i-1$, and hence \st{\tile\fka} lies weakly to the right of $\fkb$.
\end{pf}

\begin{lemma}\label{leftlowest}
Suppose $\la,\mu$ are partitions with $\la\supseteq\mu$, and $T$ is a left-\cedt{} of $\skw\la\mu$. If $\fka\in\skw\la\mu$ is the lowest node in its column, then every node in \tile\fka{} to the right of $\fka$ is the lowest node in its column.
\end{lemma}

\begin{pf}
Suppose $\fka$ lies in column $i$, and proceed by induction on the number of nodes to the right of $\fka$ in \tile\fka. Assuming $\fka$ is not the end node of its tile, there is a node $\fkb$ in \tile\fka{} in column $i+1$, which must be either \ne\fka{} or \se\fka. The only way $\fkb$ can fail to be the lowest node in column $i+1$ is if $\fkb=\ne\fka$ and $\se\fka\in\skw\la\mu$. But in this case, \se\fka{} is not attached to $\fka$ or to $\sz\fka$ (which is not a node of $\skw\la\mu$), and so \se\fka{} is the start of its tile; but this contradicts \cref{cetfae}.

So $\fkb$ is the lowest node in its column. By induction every node to the right of $\fkb$ in the same tile is the lowest node in its column, and we are done.
\end{pf}

\begin{propn}\label{ls1l}
Suppose $\la,\mu$ are partitions with $\la\supseteq\mu$. Then $\skw\la\mu$ admits at most one left-\cedt, and at most one right-\cedt. If $\skw\la\mu$ admits both a left- and a right-\cedt, then these tilings coincide.
\end{propn}

\begin{pf}
We use induction on $|\skw\la\mu|$. If $\la=\mu$ then the result is trivial, so assume $\la\supset\mu$. Let $\fka$ be the unique leftmost node of $\skw\la\mu$, and suppose $\fka$ lies in column $i$. Let $m\gs0$ be maximal such that columns $i,\dots,i+m$ each contain a node of height at most \he\fka. Suppose $T$ is a left-\cedt{} of $\skw\la\mu$. 
\clamn1
\tile\fka{} consists of the lowest nodes in columns $i,\dots,i+m$.
\prof
Since $\fka$ is the start node of its tile, every node in \tile\fka{} has height at most \he\fka, and in particular \tile\fka{} cannot contain a node in column $i+m+1$ or further to the right.

Now we prove by induction on $l$ that \tile\fka{} contains the lowest node in column $l$, for $i\ls l\ls i+m$. Let $\fkb$ be the lowest node in column $l$; then (assuming $l>i$) the lowest node in column $l-1$ is either \nw\fkb{} or \sw\fkb. In the first case, $\fkb$ cannot be the start of its tile, by \cref{cetfae}; $\fkb$ cannot be attached to \sw\fkb{}, since this is not a node of $\skw\la\mu$, and so $\fkb$ is attached to \nw\fkb, and hence lies in \tile\fka. In the second case, the fact that $\he{\sw\fkb}<\he\fkb\ls\he\fka$ means that \sw\fkb{} cannot be the end node of \tile\fka, so is attached to either $\fkb$ or $\sz\fkb$; but \sz\fkb{} is not a node of $\skw\la\mu$, and hence \sw\fkb{} is attached to $\fkb$, i.e.\ $\fkb\in\tile\fka$.
\malc
The definition of $m$ means that the nodes in \tile\fka{} can be removed to leave a smaller skew Young diagram $\skw\la\nu$, and the fact that $T$ is left-cover-expansive means that $T\setminus\{\tile\fka\}$ is a left-\cedt{} of $\skw\la\nu$. By induction on $|\skw\la\mu|$ there is at most one such tiling, and so $T$ is uniquely determined.

So there is at most one left-\cedt{} of $\skw\la\mu$, and similarly at most one right-\cedt. To prove the final statement, we continue to assume that $\skw\la\mu$ is non-empty; we choose a connected component $\calc$ of $\skw\la\mu$, and let $\fka$ denote the unique leftmost node of $\calc$, and $\fkc$ the unique rightmost node of $\calc$. (So there is no node in the column to the left of the column containing $\fka$ or in the column to the right of the column containing $\fkc$, but there are nodes in all columns in between.)
\clamn2
Suppose there exists a left-\cedt{} $T$ of $\skw\la\mu$. Then $\he\fka\ls\he\fkc$, and if equality occurs then $\fka$ and $\fkc$ lie in the same tile in $T$.
\prof
Let $\fkb_1,\dots,\fkb_r$ be the nodes in $\calc$ which are both the end nodes of their tiles and the lowest nodes in their columns, numbering them so that they appear in order from left to right. Then $\fkb_1=\en{\tile\fka}$ by Claim 1, and $\fkb_r=\fkc$. We claim that $\he{\fkb_1}<\dots<\he{\fkb_r}$. Given $1\ls l<r$, let $\fkd$ be the lowest node in the column to the right of $\fkb_l$. Then $\fkd$ is either \ne{\fkb_l} or \se{\fkb_l}; but in the latter case, $\fkd$ must be the start node of its tile, and this contradicts \cref{citfae}. So $\fkd=\ne{\fkb_l}$, and in particular $\he\fkd>\he{\fkb_l}$. By \cref{leftlowest} $\en{\tile\fkd}=\fkb_{l+1}$, and hence $\he{\fkb_{l+1}}=\he\fkd>\he{\fkb_l}$.

So we have $\he\fka=\he{\fkb_1}<\dots<\he{\fkb_r}=\he\fkc$, so $\he\fka\ls\he\fkc$, with equality if and only if $r=1$, in which case $\fkc=\fkb_1=\en{\tile\fka}$.
\malc
Now we can complete the proof. Assume there is a left-\cedt{} $T$ of $\skw\la\mu$ and a right-\cedt{} $U$. By Claim 2, we have $\he\fka\ls\he\fkc$, and symmetrically (since $U$ exists) we have $\he\fka\gs\he\fkc$. Hence $\he\fka=\he\fkc$, and so $\fka,\fkc$ lie in the same tile $t\in T$, which consists of the lowest node in every column of $\calc$. Similarly, $t$ is a tile in $U$. Removing $t$ from $\skw\la\mu$ yields a smaller skew Young diagram, and $T\setminus\{t\}$ is a left-\cedt{} of this diagram, while $U\setminus\{t\}$ is a right-\cedt. By induction $T\setminus\{t\}=U\setminus\{t\}$, and hence $T=U$.
\end{pf}

Now we restrict attention to \cedt s. We shall need the following lemma, which examines the effect of the cover-expansive property on depths of nodes.

\begin{lemma}\label{cicedepth}
Suppose $\la\supseteq\mu$, $T$ is a \cedt{} of $\skw\la\mu$, and $\fka,\nz\fka\in\skw\la\mu$. Then $\de{\nz{\fka}}<\de\fka$.
\end{lemma}

\begin{pf}
Suppose the lemma is false, and take $\fka\in\skw\la\mu$ as far to the left as possible such that $\nz\fka\in\skw\la\mu$ and $\de{\nz\fka}\gs\de\fka$. We have $\nw\fka\in\skw\la\mu$, and \nw\fka{} cannot be the end node of its tile by \cref{cetfae}, so is attached to either $\fka$ or \nz\fka.

If \nw\fka{} is attached to \nz\fka, then $\fka$ must be attached to \sw\fka{} since it cannot be the start node of its tile (again by \cref{cetfae}). But now $\de{\nw\fka}-\de{\sw\fka}=(\de{\nz\fka}+1)-(\de\fka+1)\gs0$, contradicting the choice of $\fka$.

So suppose instead that \nw\fka{} is attached to $\fka$. Then $\nz\fka$ cannot be the start node of its tile, since $\de{\nz\fka}\gs\de\fka=\de{\nw\fka}+1>0$. So \nz\fka{} is attached to \nw{\nz\fka}. But now $\de{\nw{\nz\fka}}-\de{\nw\fka}=(\de{\nz\fka}-1)-(\de\fka-1)\gs0$, and again the choice of $\fka$ is contradicted.
\end{pf}

We remark that, in contrast to the similar condition in \cref{citfae}(\ref{:de}) for \cidt s, the condition in \cref{cicedepth} does not imply the cover-expansive condition. For example, the unique Dyck tiling of $\skw{(2)}\varnothing$ satisfies this condition (trivially) but is not right-cover-expansive.

The following lemma will be useful in the next section.

\begin{lemma}\label{addrem}
Suppose $\la\supseteq\mu$, and that $\mu$ has an addable node in column $j$. If there exists a \cedt{} of $\skw\la\mu$, then $\la$ has either an addable or a removable node in column $j$.
\end{lemma}

\begin{pf}
Suppose not, and let $T$ be the \cedt{} of $\skw\la\mu$. There must be at least one node in column $j$ of $\skw\la\mu$ (otherwise the addable node of $\mu$ would also be an addable node of $\la$). Let $\fka$ be the highest node in column $j$, and consider the nodes attached to $\fka$ in $T$.

The lowest node in column $j$ of $\skw\la\mu$, namely $\fkm$, is not attached to its SW neighbour (since this is not a node of $\skw\la\mu$). Now let $\fkb$ be the highest node in column $j$ which is not attached to its SW neighbour, and suppose $\fkb\neq\fka$. Then the choice of $\fkb$ means that \nw\fkb{} is attached to \nz\fkb, so $\fkb$ is the start node of its tile; but $\nw\fkb\in\skw\la\mu$, and this contradicts \cref{cetfae}. So $\fkb=\fka$.

So $\fka$ is not attached to \sw\fka, and symmetrically is not attached to \se\fka. If $\fka$ is attached to either \nw\fka{} or \ne\fka, then (since $T$ is a Dyck tiling) it is attached to both \nw\fka{} and \ne\fka, and hence both of these nodes belong to $\la$, so $\la$ has an addable node in column $j$. The remaining possibility is that $\fka$ is a singleton tile. But now \cref{cetfae} implies that neither \nw\fka{} nor \ne\fka{} is a node of $\skw\la\mu$, so $\la$ has a removable node in column $j$.
\end{pf}

\subsection{Recurrences}

We now reassume the notation in \cref{cidtsec}: $\la,\mu$ are partitions with addable nodes $\fkl,\fkm$ respectively in column $j$, and ${\la^+},\mu^+$ are the partitions obtained by adding these nodes.

\begin{propn}\label{exp1}
Suppose $\la\supseteq\mu^+$, and let $\fka=\sz\fkl$. Then:
\begin{enumerate}
\item
there exists a \cedt{} of $\skw\la\mu$ in which $\fka$ has depth $1$ if and only if there exists a \cedt{} of $\skw\la{\mu^+}$;
\item
there exists a \cedt{} of $\skw\la\mu$ in which $\fka$ has depth greater than $1$ if and only if there exists a \cedt{} of $\skw{\la^+}{\mu^+}$.
\end{enumerate}
\end{propn}

\begin{pf}
Suppose there is a \cedt{} $T$ of $\skw\la\mu$. Arguing as in the proof of \cref{addrem}, $\fka$ cannot be attached to either \se\fka{} or \sw\fka{} (since $\mu$ has an addable node in column $j$), and $\fka$ cannot be a singleton tile (since $\la$ has an addable node in column $j$). So $\fka$ is attached to both \nw\fka{} and \ne\fka. So by the cover-expansive conditions, every node in column $j$ is attached to its NE and NW neighbours.
\begin{enumerate}
\item\label{bj1}
Suppose $\fka$ has depth $1$ in $T$. Construct a \cedt{} of $\skw\la{\mu^+}${} $U$ from $T$ as follows:
\begin{itemize}
\item
detach each node in column $j$ from its NW and NE neighbours;
\item
remove $\fkm$;
\item
attach each node in column $j$ to its SE and SW neighbours.
\end{itemize}
It is easy to check that $U$ really does give a \cedt; the tiles in $U$ but not in $T$ are obtained as follows:
\begin{itemize}
\item
if $u$ is the tile containing $\fka$ in $T$, then the portions of $u$ ending at \nw\fka{} and starting at \ne\fka{} are both tiles in $U$; since $\fka$ has depth $1$ in $T$, these are both Dyck tiles;
\item
if $u$ is a tile in $T$ containing a node $\fkc\neq\fka$ in column $j$, then $U$ contains the tile obtained by replacing $\fkc$ with \nz\fkc{} in $u$; this is a Dyck tile, since $\fkc$ has depth greater than $1$ in $T$ by \cref{cicedepth}.
\end{itemize}
So $U$ is a Dyck tiling, and the cover-expansive property follows easily from the corresponding property of $T$.

For the other direction, suppose $T$ is a \cedt{} of $\skw\la{\mu^+}$. We claim that every node in column $j$ must be attached to its SW and SE neighbours.  If there are no nodes in column $j$ (i.e.\ if $\fkm=\sz\fkl$) then this statement is trivial, so suppose otherwise; then $\fka\in\skw\la{\mu^+}$. $\fka$ cannot be attached to \nw\fka{} in $T$, since then every node in column $j-1$ would be attached to its SE neighbour; but the SE neighbour of the bottom node in column $j-1$ is not a node of $\skw\la{\mu^+}$. Since $\fka$ is not attached to \nw\fka, it must be attached to \sw\fka. Similarly $\fka$ is attached to \se\fka, and the cover-expansive property implies that every node in column $j$ is attached to its SE and SW neighbours as claimed. Now we construct a \cedt{} of $\skw\la\mu$ as follows:
\begin{itemize}
\item
detach each node in column $j$ from its SW and SE neighbours;
\item
add $\fkm$;
\item
attach each node in column $j$ to its NE and NW neighbours.
\end{itemize}
Again, it is easy to see that we have a \cedt. Moreover, $\fka$ has depth $1$ in this tiling, since \nw\fka{} is the end node of its tile in $T$, and so has depth $0$ (in both tilings).
\item\label{bj2}
Suppose $\fka$ has depth greater than $1$ in $T$. Construct a \cedt{} $U$ of $\skw{\la^+}{\mu^+}$ from $T$ as follows:
\begin{itemize}
\item
detach each node in column $j$ from its NW and NE neighbours;
\item
add $\fkl$ and remove $\fkm$;
\item
attach each node in column $j$ to its SE and SW neighbours.
\end{itemize}
Again, $U$ is a Dyck tiling, since every node in column $j$ of $\skw\la\mu$ has depth greater than $1$ in $T$. And the cover-expansive property for $U$ follows from that for $T$.

For the other direction, suppose $T$ is a \cedt{} of $\skw{\la^+}{\mu^+}$. Arguing as in (\ref{bj1}), every node in column $j$ is attached to its SW and SE neighbours. Now we construct a \cedt{} of $\skw\la\mu$ $U$ as follows:
\begin{itemize}
\item
detach each node in column $j$ from its SW and SE neighbours;
\item
add $\fkm$ and remove $\fkl$;
\item
attach each node in column $j$ to its NE and NW neighbours.
\end{itemize}
Again, it is easy to see that $U$ is a \cedt. Moreover, since \nw\fka{} is attached to $\fkl$ in $T$, it has positive depth (in both tilings), so $\fka$ has depth at least $2$ in $U$.\qedhere
\end{enumerate}
\end{pf}

\begin{cory}\label{etrec1}
$\et\la\mu=\et\la{\mu^+}+\et{\la^+}{\mu^+}$.
\end{cory}

\begin{pf}
First suppose $\la\supseteq\mu^+$. The first paragraph of the proof of \cref{exp1} shows that the node $\fka=\sz\fkl$ must have positive depth in a \cedt{} of $\skw\la\mu$. Hence the result follows from \cref{exp1}.

Alternatively, suppose $\la\nsupseteq\mu^+$. Then we have $\la\supseteq\mu$ if and only if $\la^+\supseteq\mu^+$, and if these conditions hold then $\skw\la\mu=\skw{\la^+}{\mu^+}$; so $\et\la\mu=\et{\la^+}{\mu^+}$.
\end{pf}

\begin{propn}\label{exp2}
$\et{\la^+}\mu=\et\la\mu$.
\end{propn}

\begin{pf}
Suppose $T$ is a \cedt{} of $\skw{\la^+}\mu$. Then $\fkl$ cannot be attached to \se\fkl, since then every node in column $j$ would be attached to its SE neighbour; but \se\fkm{} is not a node of $\skw{\la^+}\mu$. Similarly $\fkl$ is not attached to \sw\fkl, so forms a singleton tile in $T$. Removing this tile yields a \cedt{} of $\skw\la\mu$.

In the other direction, suppose $T$ is a \cedt{} of $\skw\la\mu$. From the proof of \cref{exp1}, $\sz\fkl$ is attached to both \se\fkl{} and \sw\fkl. Hence if we add $\fkl$ as a singleton tile, the resulting tiling is a \cedt{} of $\skw{\la^+}\mu$.
\end{pf}

\section{Young permutation modules for two-part compositions}\label{ypmsec}

In this section we apply our results on Dyck tilings to compute change-of-basis matrices for a certain module for the symmetric group.

\subsection{The Young permutation module $\mfg$}\label{ypmsubsec}

Suppose $k\gs0$, and let $\sss k$ denote the symmetric group of degree $k$, with $\{t_1,\dots,t_{k-1}\}$ the set of Coxeter generators (so $t_i$ is the transposition $\perc(i\ipo)$). Let $\bbf$ be any field, and consider the group algebra $\bbf\sss k$. Define $s_i=t_i-1\in\bbf\sss k$ for $i=1,\dots,n-1$.

Now write $k=f+g$ with $f,g$ non-negative integers. Consider the \emph{Young permutation module} $\mfg$ for $\bbf\sss k$ indexed by the composition $(f,g)$. This is just the permutation module on the set of cosets of the maximal Young subgroup $\sss f\times\sss g$, and has the following presentation:
\[
\mfg=\left\lan\gn\ \left|\ t_i\gn=\gn\text{ for }i\neq f\right.\right\ran.
\]

\begin{lemma}\label{ann}
For $1\ls i\ls k-2$, the element $s_is_{i+1}s_i-s_i=t_is_{i+1}s_i-s_{i+1}s_i-s_i$ annihilates $\mfg$.
\end{lemma}

\begin{pf}
In terms of permutations, the given element is
\[
\perc(i\ipt)-\perc(i\ipo\ipt)-\perc(i\ipt\ipo)+\perc(i\ipo)+\perc(\ipo\ipt)-1.
\]
$\mfg$ may be viewed as the permutation module on the set of $f$-subsets of the set $\{1,\dots,k\}$; it is easily seen that the given element kills any such subset.
\end{pf}

$\mfg$ has a basis indexed by the set of minimal left coset representatives of $\sss f\times\sss g$ in $\sss k$; this set is in one-to-one correspondence with the set $\lfg$ of partitions $\la$ for which $\la_1\ls f$ and $\la_1'\ls g$; given $\la\in\lfg$, we write the corresponding basis element as $t_\la\gn$, where
\begin{align*}
t_\la&=\left(t_{\la_1+g}t_{\la_1+g+1}\dots t_{f+g-1}\right)\left(t_{\la_2+g-1}t_{\la_2+g}\dots t_{f+g-2}\right)\dots\left(t_{\la_g+1}t_{\la_g+2}\dots t_f\right).
\\
\intertext{Our objective here is to study the elements}
s_\la&=\left(s_{\la_1+g}s_{\la_1+g+1}\dots s_{f+g-1}\right)\left(s_{\la_2+g-1}s_{\la_2+g}\dots s_{f+g-2}\right)\dots\left(s_{\la_g+1}s_{\la_g+2}\dots s_f\right)
\end{align*}
for $\la\in\lfg$. It is easy to see that $\rset{s_\la\gn}{\la\in\lfg}$ is also a basis for $\mfg$, since when each each $s_\la\gn$ is expressed as a linear combination of the elements $t_\mu\gn$, the matrix of coefficients is unitriangular with respect to a suitable ordering.

We shall use cover-expansive Dyck tilings to describe this change-of-basis matrix explicitly, and then describe its inverse using cover-inclusive Dyck tilings. We also give a simple expression for the sum $\sum_{\la\in\lfg}t_\la\gn$ as a linear combination of the elements $s_\la\gn$, which will be useful in \cref{klrsec}.

\subsection{Change of basis}

Our first result on change-of-basis matrices is the following.

\begin{thm}\label{mainsn1}
Suppose $f,g$ are non-negative integers, and $\mu\in\lfg$. Then
\[
s_\mu\gn=\displaystyle\zsum{\la\in\lfg}(-1)^{|\la|+|\mu|}\et\la\mu t_\la\gn.
\]
\end{thm}

We begin with some simple observations concerning the actions of the generators $t_1,\dots,t_{k-1}$ on the basis elements $t_\la$. We continue to use the Russian convention for Young diagrams.

\begin{lemma}\label{actiontt}
Suppose $\la\in\lfg$ and $1\ls i<k$. Then:
\begin{itemize}
\item
if $\la$ has an addable node $\fkl$ in column $i-g$, then $t_it_\la\gn=t_{\la^+}\gn$, where $\la^+$ is the partition obtained by adding $\fkl$ to $\la$;
\item
if $\la$ has a removable node $\fkl$ in column $i-g$, then $t_it_\la\gn=t_{\la^-}\gn$, where $\la^-$ is the partition obtained by removing $\fkl$ from $\la$;
\item
if $\la$ has neither an addable nor a removable node in column $i-g$, then $t_it_\la\gn=t_\la\gn$.
\end{itemize}
\end{lemma}

\begin{pf}
This is an easy consequence of the definitions and the Coxeter relations.
\end{pf}

\begin{pf}[Proof of \cref{mainsn1}]
We proceed by downwards induction on $|\mu|$. If $\mu=(f^g)$, then $s_\mu=t_\mu=1$ and the result follows. Assuming $\mu\neq(f^g)$, $\mu$ has an addable node in column $j$, for some $-g<j<f$. We let ${\mu^+}$ denote the partition obtained by adding this addable node; then we have $s_\mu=s_{j+g}s_{\mu^+}$, and by induction
\[
s_{\mu^+}\gn=\zsum{\la\in\lfg}(-1)^{|\la|+|{\mu^+}|}\et\la{\mu^+}t_\la\gn.
\]
Let $\lfgp$ denote the set of $\la\in\lfg$ having an addable node in column $j$, and for $\la\in\lfgp$ let $\la^+$ denote the partition obtained by adding this addable node; similarly, let $\lfgm$ denote the set of partitions in $\lfg$ having a removable node in column $j$, and for each $\la\in\lfgm$ let $\la^-$ denote the partition obtained by removing this removable node. Note that the functions $\la\mapsto\la^+$ and $\la\mapsto\la^-$ define mutually inverse bijections between $\lfgp$ and $\lfgm$. Now
\begin{align*}
s_\mu\gn&=(t_{j+g}-1)\zsum{\la\in\lfg}(-1)^{|\la|+|{\mu^+}|}\et\la{\mu^+} t_\la\gn\\
&=\zsum{\la\in\lfgp}(-1)^{|\la|+|{\mu^+}|}\et\la{\mu^+}(t_{\la^+}-t_\la)\gn+\zsum{\la\in\lfgm}(-1)^{|\la|+|{\mu^+}|}\et\la{\mu^+}(t_{\la^-}-t_\la)\gn\tag*{by \cref{actiontt}}\\
&=\zsum{\la\in\lfgp}(-1)^{|\la|+|\mu|}\big(\et\la{\mu^+}+\et{\la^+}{\mu^+}\big)t_\la\gn+\zsum{\la\in\lfgm}(-1)^{|\la|+|\mu|}\big(\et\la{\mu^+}+\et{\la^-}{\mu^+}\big)t_\la\gn\\
&=\zsum{\la\in\lfgp}(-1)^{|\la|+|\mu|}\et\la\mu t_\la\gn+\zsum{\la\in\lfgm}(-1)^{|\la|+|\mu|}\et{\la^-}\mu t_\la\gn\tag*{by \cref{exp1}}\\
&=\zsum{\la\in\lfgp\cup\lfgm}(-1)^{|\la|+|\mu|}\et\la\mu t_\la\gn\tag*{by \cref{exp2}}\\
&=\zsum{\la\in\lfg}(-1)^{|\la|+|\mu|}\et\la\mu t_\la\gn\tag*{by \cref{addrem}.}
\end{align*}
So \cref{mainsn1} follows by induction.
\end{pf}

Now we give our second main result on change-of-basis matrices.

\begin{thm}\label{mainsn2}
Suppose $f,g$ are non-negative integers, and $\mu\in\lfg$. Then
\[
t_\mu\gn=\displaystyle\zsum{\la\in\lfg}\ndt\la\mu s_\la\gn.
\]
\end{thm}

The proof of this result is rather more difficult. To begin with, we compute the actions of $s_1,\dots,s_{k-1}$ on the basis elements $s_\la\gn$.

\begin{propn}\label{actionts}
Suppose $\mu\in\lfg$ and $1\ls i<k$. Then exactly one of the following occurs.
\begin{enumerate}
\item
$\mu$ has an addable node in column $i-g$. In this case $s_is_\mu\gn=-2s_\mu\gn$.
\item
$\mu$ has a removable node in column $i-g$. In this case $s_is_\mu\gn=s_{\mu^-}\gn$, where $\mu^-$ denotes the partition obtained by removing this node.
\item
For some $0\ls a\ls g$ we have $\mu_a>i-g+a>\mu_{a+1}$ (where the left-hand inequality is regarded as automatically true in the case $a=0$).
\begin{enumerate}
\item
If $\mu_w<i-g+2a-w$ for all $a<w\ls g$, then $s_is_\mu\gn=0$.
\item
Otherwise, let $w>a$ be minimal such that $\mu_w=i-g+2a-w$, and set
\[
\mu^{a,w}=(\mu_1,\dots,\mu_a,i-g+a,\mu_{a+1}+1,\dots,\mu_{w-1}+1,\mu_{w+1},\dots,\mu_g).
\]
Then $s_is_\mu\gn=s_{\mu^{a,w}}\gn$.
\end{enumerate}
\item
For some $1\ls a<g$ we have $\mu_a=\mu_{a+1}=i-g+a$.
\begin{enumerate}
\item
If $i+2a>k$ and $\mu_w<i-g+2a-w$ for $w=1,\dots,a-1$, then $s_is_\mu\gn=0$.
\item
Otherwise, let $w<a$ be maximal such that $\mu_w\gs i-g+2a-w$ (taking $w=0$ if there is no such $w$), and define
\[
\mu^{w,a}=(\mu_1,\dots,\mu_w,i-g+2a-w,\mu_{w+1}+1,\dots,\mu_{a-1}+1,\mu_{a+1},\dots,\mu_g).
\]
Then $s_is_\mu\gn=s_{\mu^{w,a}}\gn$.
\end{enumerate}
\end{enumerate}
\end{propn}

\begin{egs}\indent
\begin{enumerate}
\vspace{-\topsep}
\item
Suppose $f=4$, $g=5$ and $\mu=(4,2^4)$. Then
\[
s_\mu=s_6s_7s_5s_6s_4s_5s_3s_4.
\]
Taking $i=3$, we find that $\mu$ satisfies condition 4(a) of \cref{actionts}, with $a=4$. And indeed
\begin{align*}
s_3s_\mu m&=s_3s_6s_7s_5s_6s_4s_5s_3s_4m\\
&=s_6s_7s_5s_6(s_3s_4s_3)s_5s_4m\\
&=s_6s_7s_5s_6s_3s_5s_4m\tag*{by \cref{ann}}\\
&=s_6s_7(s_5s_6s_5)s_3s_4m\\
&=s_6s_7s_5s_3s_4m\tag*{by \cref{ann}}\\
&=s_6s_5s_3s_4s_7m\\
&=0\tag*{since $s_7m=0$.}
\end{align*}
\item
Now suppose $f=4$, $g=5$ and $\mu=(4,2,1^3)$. Then
\[
s_\mu=s_6s_7s_4s_5s_6s_3s_4s_5s_2s_3s_4.
\]
Taking $i=2$, we find that $\mu$ satisfies condition 4(b) of \cref{actionts}, with $a=4$ and $w=1$, giving $\mu^{w,a}=(4^2,3,2,1)$. Now we have
\begin{align*}
s_2s_\mu m&=s_2s_6s_7s_4s_5s_6s_3s_4s_5s_2s_3s_4m\\
&=s_6s_7s_4s_5s_6(s_2s_3s_2)s_4s_5s_3s_4m\\
&=s_6s_7s_4s_5s_6s_2s_4s_5s_3s_4m\\
&=s_6s_7(s_4s_5s_4)s_6s_2s_5s_3s_4m\\
&=s_6s_7s_4s_6s_2s_5s_3s_4m\\
&=(s_6s_7s_6)s_4s_2s_5s_3s_4m\\
&=s_6s_4s_2s_5s_3s_4m\\
&=m_{(4^2,3,2,1)}m.
\end{align*}
\end{enumerate}
\end{egs}

\begin{pf}[Proof of \cref{actionts}]
To see that $\mu$ satisfies exactly one of the four conditions, let $\fkl$ denote the lowest node in column $j$ which is not a node of $\mu$, and consider whether \se\fkl{} and \sw\fkl{} lie in $\mu$; for the purposes of this argument we regard all nodes of the form $(a,0)$ or $(0,b)$ as lying in $\mu$. If both \se\fkl{} and \sw\fkl{} lie in $\mu$, then $\fkl$ is an addable node and we are in case (1). If neither lies in $\mu$, then \sz\fkl{} is a removable node, and we are in case (2). If $\se\fkl\in\mu\notni\sw\fkl$, then we are in case (3), while if $\se\fkl\notin\mu\ni\sw\fkl$ then we are in case (4).

Now we analyse the four cases.
\begin{enumerate}
\item
Since $\mu$ has an addable node in column $i-g$, $s_\mu$ can be written in the form $s_is_{\mu^+}$. We have $s_i^2=-2s_i$, and the result follows.
\item
By definition (and the fact that $s_i$ and $s_j$ commute when $j\neq i\pm1$) $s_{\mu^-}=s_is_\mu$.
\item\label{case3}
\begin{enumerate}
\item
We proceed by induction on $i$. Consider the case $a=g$. In this case $s_\mu$ only involves terms $s_j$ for $j\gs i+2$, so (since $s_i\gn=0$) we get $s_is_\mu\gn=0$.

Now suppose $a<g$, and set
\[
\bar\mu=(\mu_1,\dots,\mu_a,i-g+a,\mu_{a+2},\dots,\mu_g).
\]
Then $\bar\mu$ satisfies the inductive hypothesis, with $a$ replaced by $a+1$ and $i$ replaced by $i-2$. Furthermore, we can write
\[
s_\mu=s_{\mu_{a+1}+g-a}\dots s_{i-1}s_{\bar\mu}.
\]
$\bar\mu$ has an addable node $(a+1,i-g+a+1)$ in column $i-g$, so $s_{\bar\mu}\gn$ can be written in the form $s_ix$ for some $x\in\mfg$. So we have
\begin{align*}
s_is_\mu\gn&=s_is_{\mu_{a+1}+g-a}\dots s_{i-1}s_ix\\
&=s_{\mu_{a+1}+g-a}\dots s_{i-2}s_is_{i-1}s_ix\\
&=s_{\mu_{a+1}+g-a}\dots s_{i-2}s_ix\tag*{by \cref{ann}}\\
&=s_{\mu_{a+1}+g-a}\dots s_{i-2}s_{\bar\mu}\gn\\
&=0\tag*{by induction.}
\end{align*}
\item
We use induction on $w-a$. If $w=a+1$, then we have $s_\mu=s_{i-1}s_{\mu^{a,w}}$, and $\mu^{a,w}$ has an addable node $(a+1,i-g+a+1)$ in column $i-g$, so $s_{\mu^{a,w}}\gn$ can be written as $s_ix$ for some $x\in\mfg$. So by \cref{ann}
\[
s_is_\mu\gn=s_is_{i+1}s_ix=s_ix=s_{\mu^{a,w}}\gn.
\]

Now assume that $w>a+1$, and as above set
\[
\bar\mu=(\mu_1,\dots,\mu_a,i-g+a,\mu_{a+2},\dots,\mu_g).
\]
Then $\bar\mu$ satisfies the inductive hypothesis, with $a$ replaced by $a+1$ and $i$ replaced by $i-2$, and with the same value of $w$, yielding
\[
\bar\mu^{a+1,w}=(\mu_1,\dots,\mu_a,i-g+a,i-g+a-1,\mu_{a+2}+1,\dots,\mu_{w-1}+1,\mu_{w+1},\dots,\mu_g).
\]
As above we write $s_{\bar\mu}=s_ix$, and obtain
\begin{align*}
s_is_\mu\gn&=s_{\mu_{a+1}+g-a}\dots s_{i-2}s_{\bar\mu}\gn\\
&=s_{\mu_{a+1}+g-a}\dots s_{i-3}s_{\bar\mu^{a+1,w}}\gn\tag*{by induction}\\
&=s_{\mu^{a,w}}\gn.
\end{align*}
\end{enumerate}
\item
This is symmetrical to case (\ref{case3}), replacing partitions with their conjugates and swapping $f$ and $g$.\qedhere
\end{enumerate}
\end{pf}

We now seek to re-phrase \cref{actionts} so that for a given $\la,\mu$ we can write down the coefficient of $s_\la\gn$ in $s_is_\mu\gn$. Fix $i$, and suppose that as in part (3) of \cref{actionts} we have $\mu_a>i-g+a>\mu_{a+1}$ for some $1\ls a\ls g$. Let $w>a$ be minimal such that $\mu_w=i-g+2a-w$, assuming there is such a $w$, and write $\mu\rea i\mu^{a,w}$, where $\mu^{a,w}$ is as defined in \cref{actionts}.

Now recall some notation from \cref{defnsec}: if $j$ is fixed and $\la$ is a partition with an addable node $\fkl$ in column $j$, then $\xx\la$ denotes the set of all integers $x$ such that $\la$ has a removable node $\fkn$ in column $j+x$, with $\he\fkn=\he\fkl-1$, and $\he\fkp<\fkl$ for all nodes $\fkp$ in all columns between $j$ and $j+x$. $\xxp\la$ denotes the set of positive elements of $\xx\la$, and $\xxm\la$ the set of negative elements. Given $x\in\xxp\la$, $\rrh\la x$ is the partition obtained from $\la$ by removing the highest node in each column from $j+1$ to $j+x$; $\rrh\la x$ is defined similarly for $x\in\xxm\la$.

\begin{lemma}\label{harp}
Suppose $\la,\mu\in\lfg$ and $1\ls i<k$, and set $j=i-g$. 
Then $\mu\rea i\la$ if and only if $\la$ has an addable node in column $j$ and $\mu=\rrh\la x$ for some $x\in\xxm\la$.
\end{lemma}

\begin{pf}\indent
\begin{description}
\vspace{-\topsep}
\item[$(\Rightarrow)$]
If $\mu\rea i\la$ then we have
\[
\la=\mu^{a,w}=(\mu_1,\dots,\mu_a,j+a,\mu_{a+1}+1,\dots,\mu_{w-1}+1,\mu_{w+1},\dots,\mu_g),
\]
where $\mu_a>j+a>\mu_{a+1}$ and $w>a$ is minimal such that $\mu_w\gs j+2a-w$. Observe that $\la$ has an addable node $\fkl=(a+1,j+a+1)$ in column $j$. The choice of $w$ implies that $\mu_{w-1}=\mu_w=j+2a-w$, and so $\la$ has a removable node $\fkn=(w,j+2a-w+1)$ in column $j+x$, where $x=2(a-w)+1$. Furthermore, $\he\fkn=j+2a=\he\fkl-1$, and we claim that all the nodes in columns $j+x+1,\dots,j-1$ of $\la$ have height at most $j+2a$: if this condition fails, then it must fail for a node of the form $(x,\la_x)$ with $a+2\ls x\ls w-1$, but for $x$ in this range we have $\he{(x,\la_x)}=x+\la_x=x+\mu_{x-1}+1\ls x+j+2a-x$. So $x\in\xxm\la$.

To construct $\mu$ from $\la$ we remove the nodes
\begin{align*}
&(w,\mu_w+1),(w,\mu_w+2),\dots,(w,\mu_{w-1}+1),\\
&(w-1,\mu_{w-1}+1),(w-1,\mu_{w-1}+2),\dots,(w-1,\mu_{w-2}+1),\\
&\dots\\
&(a+2,\mu_{a+2}+1),(a+2,\mu_{a+2}+2),\dots,(a+2,\mu_{a+1}+1),\\
&(a+1,\mu_{a+1}+1),(a+1,\mu_{a+1}+2),\dots,(a+1,j+a),
\end{align*}
which lie in columns $j+x,\dots,j-1$ respectively. So $\mu=\rrh\la x$.
\item[$(\Leftarrow)$]
Suppose $\la$ has an addable node in column $j$. Writing this node as $(a+1,\la_{a+1}+1)$, we have $a=0$ or $\la_a>\la_{a+1}$, and $\la_{a+1}-a=j$.

Now suppose $x\in\xxm\la$. Then $\la$ has a removable node $\fkn=(w,\la_w)$ in column $j+x$, and $w>a$ since $x<0$. Then $\la_w>\la_{w+1}$ and $\la_w-w=j+x$, and (since $\he\fkn=\he\fkl-1$) $w+\la_w=a+1+\la_{a+1}$. The fact that $\he\fkp<\he\fkl$ for all nodes $\fkp$ in columns $j+x+1,\dots,j-1$ implies in particular that $\he{(x,\la_x)}<\he\fkl$ for $a+2\ls x\ls w-1$, i.e.\ $x+\la_x\ls a+1+\la_{a+1}$.

Constructing $\rrh\la x$ involves removing the nodes
\begin{align*}
(w,\la_w),{}&(w-1,\la_w),(w-1,\la_w+1),\dots,(w-1,\la_{w-1}),\\
&(w-2,\la_{w-1}),(w-2,\la_{w-1}+1),\dots,(w-2,\la_{w-2}),\\
&\dots\\
&(a+1,\la_{a+2}),(a+1,\la_{a+2}+1),\dots,(a+1,\la_{a+1}),
\end{align*}
so
\[
\rrh\la x=(\la_1,\dots,\la_a,\la_{a+2}-1,\dots,\la_w-1,\la_w-1,\la_{w+1},\dots,\la_g).
\]
Setting $\mu=\rrh\la x$, the (in)equalities observed above for $\la$ give $\mu_a>j+a>\mu_{a+1}$, $\mu_x<j+2a-x$ for $a<x<w$ and $\mu_w=j+2a-w$. Furthermore, $\la$ equals the partition $\mu^{a,w}$, so $\mu\rea i\la$.
\qedhere
\end{description}
\end{pf}

We now note a counterpart to this for part (4) of \cref{actionts}, which follows by conjugating partitions. Suppose that $\mu\in\lfg$ satisfies the conditions of (4), and that $\mu^{w,a}$ is defined, and write $\mu\aer i\mu^{w,a}$.

\begin{lemma}\label{prah}
Suppose $\la,\mu\in\lfg$ and $1\ls i<k$, and set $j=i-g$. Then $\mu\aer i\la$ if and only if $\la$ has a removable node in column $j$ and $\mu=\rrh\la x$ for some $x\in\xxp\la$.
\end{lemma}

The following result is now immediate from \cref{actionts,harp,prah}.

\begin{cory}\label{appear}
Suppose $\mu\in\lfg$, and write $s_is_\mu\gn$ as a linear combination $\zsum{\la\in\lfg}a_\la s_\la\gn$. Then $a_\la$ equals:
\begin{itemize}
\item
$-2$, if $\la$ has an addable node in column $i-g$ and $\mu=\la$;
\item
$1$, if $\la$ has an addable node in column $i-g$, and
\begin{itemize}
\item
$\mu$ is the partition obtained by adding this node, or
\item
$\mu=\rrh\la x$ for some $x\in\xx\la$;
\end{itemize}
\item
$0$, otherwise.
\end{itemize}
\end{cory}

\begin{pf}[Proof of \cref{mainsn2}]
As with \cref{mainsn1}, we proceed by downwards induction on $|\mu|$. If $\mu=(f^g)$, then $s_\mu=t_\mu$ and the result follows, since clearly $\ndt\la\la=1$.

Assuming $\mu\neq(f^g)$, $\mu$ has an addable node in column $j$ for some $-g<j<f$. We let $\mu^+$ denote the partition obtained by adding this addable node. Then $t_\mu=t_{j+g}t_{\mu^+}$, and by induction
\[
t_{\mu^+}\gn=\zsum{\la\in\lfg}\ndt\la{\mu^+} s_\la\gn.
\]
As before, we write $\lfgp$ for the set of $\la\in\lfg$ having an addable node in column $j$, and for $\la\in\lfgp$ we let $\la^+$ denote the partition obtained by adding this addable node. We have
\begin{align*}
t_\mu\gn&=\zsum{\la\in\lfg}\ndt\la{\mu^+}(s_{j+g}+1)s_\la\gn\\
&=\zsum{\la\in\lfg}\ndt\la{\mu^+} s_\la\gn+\zsum{\la\in\lfgp}\left(\ndt{\la^+}{\mu^+}-2\ndt\la{\mu^+}+\sum_{x\in\xx\la}\ndt{\rrh\la x}{\mu^+}\right)s_\la\gn\tag*{by \cref{appear}}\\
&=\zsum{\la\in\lfg\setminus\lfgp}\ndt\la{\mu^+} s_\la\gn+\zsum{\la\in\lfgp}\left(\ndt{\la^+}{\mu^+}-\ndt\la{\mu^+}+\sum_{x\in\xx\la}\ndt{\rrh\la x}{\mu^+}\right)s_\la\gn\\
&=\zsum{\la\in\lfg\setminus\lfgp}\ndt\la\mu s_\la\gn+\zsum{\la\in\lfgp}\ndt\la\mu s_\la\gn\tag*{by \cref{noadd,recur}}\\
&=\zsum{\la\in\lfg}\ndt\la\mu s_\la\gn.
\end{align*}
\cref{mainsn2} follows by induction.
\end{pf}

\subsection{The sum of the elements $t_\mu\gn$}

It will be critical in the next section to be able to express $\zsum{\la\in\lfg}t_\la\gn$ in terms of the basis elements $s_\la\gn$. To enable us to do this, we define an integer-valued function $\fp$ on partitions.

If $\la$ is a partition, we let $\rrs\la a$ denote the partition obtained by removing the first $a$ parts of $\la$, i.e.\ the partition $(\la_{a+1},\la_{a+2},\dots)$, and let $\ccs\la b$ denote the partition obtained by reducing all the parts by $b$ (and deleting all negative parts), i.e.\ the partition $(\max\{\la_1-b,0\},\max\{\la_2-b,0\},\dots)$. We write $\rcs\la ab=\rrs{\ccs\la b}a$.

Now we can define $\fp$ recursively. Set $\fp(\varnothing)=1$. Given a partition $\la\neq\varnothing$, take a node $(a,b)$ of $\la$ for which $a+b$ is maximal (call this a \emph{highest} node of $\la$) and define
\[
\sigma=\rrs\la a,\qquad\tau=\ccs\la b.
\]
Now set
\[
\fp(\la)=\binom{a+b}a\fp(\sigma)\fp(\tau).
\]

\begin{lemma}\label{fwelld}
$\fp(\la)$ is well-defined, i.e.\ does not depend on the choice of $(a,b)$.
\end{lemma}

\begin{pf}
Suppose $(a,b)$ and $(\hi,\hj)$ are both highest nodes of $\la$; in particular, $a+b=\hi+\hj$. Assume without loss that $a<\hi$. Let $\sigma=\rrs\la a$ and $\tau=\ccs\la b$ as above, and define $\hat\sigma,\hat\tau$ correspondingly. By induction $\fp(\sigma),\fp(\tau),\fp(\hat\sigma),\fp(\hat\tau)$ are well-defined, and we must show that
\[
\fp(\sigma)\fp(\tau)\binom{a+b}a=\fp(\hat\sigma)\fp(\hat\tau)\binom{\hi+\hj}\hi.
\]
Now $(\hi-a,\hj)$ is a highest node of $\sigma$, so we have
\[
\fp(\sigma)=\fp(\hat\sigma)\fp(\xi)\binom{\hi-a+\hj}{\hi-a},
\]
where $\xi=\rcs\la a\hj$. Also, $(a,b-\hj)$ is a highest node of $\hat\tau$, so that
\[
\fp(\hat\tau)=\fp(\xi)\fp(\tau)\binom{a+b-\hj}a.
\]
So (since $a+b=\hi+\hj$) we obtain
\[
\fp(\sigma)\fp(\tau)\fp(\xi)\binom\hi a=\fp(\hat\sigma)\fp(\hat\tau)\fp(\xi)\binom b\hj.
\]
Dividing both sides by $\fp(\xi)\hi!b!$ (which is obviously positive) and multiplying by $(a+b)!(\hi-a)!=(\hi+\hj)!(\hj-b)!$ gives the result.
\end{pf}

\begin{rmk}
In the notation of the above lemma, we have
\[
\fp(\la)=\fp(\hat\sigma)\fp(\xi)\fp(\tau)\trinom a{\hi-a}\hj,
\]
where we use $\trinom xyz$ to denote the trinomial coefficient $\mfrac{(x+y+z)!}{x!y!z!}$. We will refer to this as `factorising using the nodes $(a,b)$ and $(\hi,\hj)$'. This readily generalises to factorising using any number of highest nodes, with a corresponding multinomial coefficient.
\end{rmk}

Our main objective is to prove the following statement.

\begin{thm}\label{snsum}
Suppose $f,g$ are non-negative integers. Then
\[
\sum_{\mu\in\lfg}t_\mu\gn=\zsum{\la\in\lfg}\fp(\la)s_\la\gn.
\]
\end{thm}

We prove this using the following recurrence, in which we again use the notation $\xx\la$ and $\rrh\la x$ from \cref{defnsec}.

\begin{propn}\label{fprecur}
Assume that $\la$ is a partition with an addable node $\fkl$ in column $j$, and let $\la^+$ denote the partition obtained by adding this node. Then
\[
\fp(\la^+)+\sum_{x\in\xx\la}\fp(\rrh\la x)=2\fp(\la).
\]
\end{propn}

\begin{pf}
Write $\fkl$ as $(a,b)$. Suppose first of all that $\fkl$ is not the unique highest node of $\la^+$. Take a highest node $\fkm=(\hat a,\hat b)\neq\fkl$, and suppose without loss that $\hat a<a$. Since the height of $\fkm$ is at least the height of $\fkl$, $\hat b-\hat a-j$ is strictly larger than any element of $\xx\la$. So $\fkm$ is a highest node of $\rrh\la x$ for every $x\in\xx\la$, as well as of $\la$ and $\la^+$. By induction the result holds for the partition $(\la_{\hat a+1},\la_{\hat a+2},\dots)$, and so it holds for $\la$, factorising $\fp(\la)$ (and also $\fp(\la^+)$ and $\fp(\rrh\la x)$ for each $x$) using the node $(\hat a,\hat b)$.

So we can assume that $\fkl$ is the unique highest node of $\la^+$. Let $(a_1,b_1),\dots,(a_r,b_r)$ be the nodes of $\la$ to the right of $\fkl$ for which $a_x+b_x=a+b-1$; order these nodes so that $a_1>\dots>a_r$. Note that if $r\gs1$, then automatically $(a_1,b_1)=(a-1,b)$. Similarly, let $(c_1,d_1),\dots,(c_s,d_s)$ be the nodes to the left of $\fkl$ for which $c_x+d_x=a+b-1$.

We compute $\fp(\la^+)$ by first factorising using the node $\fkl$. We may also factorise $\fp(\la)$ using the nodes $(a-1,b)$ and $(a,b-1)$ (omitting either or both of these, if $a$ or $b$ equals $1$), and we readily obtain
\[
\fp(\la^+)=\frac{a+b}{ab}\fp(\la).\tag*{\cne}
\]
Now we consider the partitions $\rrh\la x$, for $x\in\xxp\la$. From our assumptions so far, $\xxp\la$ is the set of integers $x_i=b_i-a_i-b+a$, for $i=1,\dots,r$. Choose such an $i$, and set $\kappa=\rrh\la{x_i}$. We wish to compare $\fp(\kappa)$ with $\fp(\la)$. To do this, we define
\[
\gamma=\rcs\la{a_{i+1}}{b-1},\qquad
\delta=\rcs\kappa{a_{i+1}}{b-1}
\]
(where we put $a_{i+1}=0$ if $i=r$), and we claim that $\fp(\kappa)/\fp(\la)=\fp(\delta)/\fp(\gamma)$. If $s=0$ and $i=r$, then this is trivial since $\gamma=\la$ and $\delta=\kappa$; assuming instead that $s\gs1$ or $i<r$, factorising both $\fp(\kappa)$ and $\fp(\la)$ at the nodes $(c_1,d_1),\dots,(c_s,d_s),(a_{i+1},b_{i+1}),\dots,(a_r,b_r)$ (i.e.\ at all the highest nodes of $\kappa$) gives the result.

Now $\gamma$ has highest nodes $(a-1-a_{i+1},1)$ and $(a_i-a_{i+1},b_i-b+1)$ (which coincide if $i=1$), and we factorise $\fp(\gamma)$ at these nodes, obtaining
\[
\fp(\gamma)=\fp(\mu)\fp(\nu)\frac{(a-a_{i+1})!}{(a-1-a_i)!(a_i-a_{i+1})!},
\]
where $\mu=\rcs\gamma{a_i-a_{i+1}}1$ and $\nu=\ccs\gamma{b_i-b+1}$.
\clam
$\dfrac{\fp(\delta)}{\fp(\gamma)}=\dfrac{a_i-a_{i+1}}{(a-a_i)(a-a_{i+1})}$.
\prof
In the case where $a_i=a_{i+1}+1$, we have $\mu=\delta$ and $\nu=\varnothing$, so the factorisation of $\fp(\gamma)$ above gives the result. If $a_i>a_{i+1}+1$, then $\delta$ has a highest node $(a_i-a_{i+1}-1,b_i-b+1)$, and factorising at this node we get
\[
\fp(\delta)=\fp(\mu)\fp(\nu)\binom{a_i-a_{i+1}+b_i-b}{a_i-a_{i+1}-1}.
\]
The claim follows from this factorisation and the factorisation of $\fp(\gamma)$ above, using the fact that $a_i+b_i=a+b-1$.
\malc
This claim yields
\[
\frac{\fp(\kappa)}{\fp(\la)}=\frac{\fp(\delta)}{\fp(\gamma)}=\frac1{a-a_i}-\frac1{a-a_{i+1}},
\]
and summing over $i$ we obtain
\[
\frac{\sum_{x\in\xxp\la}\fp(\rrh\la x)}{\fp(\la)}=\frac1{a-a_1}-\frac1{a-a_{r+1}}=1-\frac1a.
\]
A symmetrical calculation applies to partitions $\rrh\la x$ for $x\in\xxm\la$, yielding
\[
\frac{\sum_{x\in\xxm\la}\fp(\rrh\la x)}{\fp(\la)}=1-\frac1b.
\]
Adding these two results together and combining with \cne{} gives the result.
\end{pf}

\begin{cory}\label{sikills}
For $1\ls i<k$ we have
\[
s_i\zsum{\la\in\lfg}\fp(\la)s_\la\gn=0.
\]
\end{cory}

\begin{pf}
As above, write $\lfgp$ for the set of $\la\in\lfg$ having an addable node in column $i-g$, and for each such $\la$ let $\la^+$ denote the partition obtained by adding this addable node. Then by \cref{appear} we have
\[
s_i\zsum{\la\in\lfg}\fp(\la)s_\la\gn=\zsum{\la\in\lfgp}\left(-2\fp(\la)+\fp(\la^+)-\sum_{x\in\xx\la}\fp(\rrh\la x)\right)s_\la\gn,
\]
and each summand on the right-hand side is zero by \cref{fprecur}.
\end{pf}

\begin{pf}[Proof of \cref{snsum}]
\Cref{sikills} shows that $\sum_\la\fp(\la)s_\la\gn$ is killed by all of $s_1,\dots,s_{k-1}$, and hence lies in a trivial submodule of $\mfg$. But $\mfg$ is a permutation module for a transitive action of $\sss k$, and so has a unique trivial submodule, spanned by the sum of the elements of the permutation basis, i.e.\ $\sum_\la t_\la\gn$. So the two sides agree up to multiplication by a scalar. On the other hand, it is clear that when the right-hand side is expressed as a linear combination of the $t_\la\gn$, the coefficient of $t_\varnothing\gn$ is $\fp(\varnothing)=1$, so in fact the two sides are equal.
\end{pf}

\subsection{Results of Kenyon, Kim, M\'esz\'aros, Panova and Wilson}

As we have mentioned, cover-inclusive Dyck tilings were introduced by Kenyon and Wilson in \cite{kw}, and in fact the main result in their paper has a bearing on the results in the present paper.

First we explain how \cedt s appear in disguise in \cite{kw}. Given a partition $\la$, we define a sequence of parentheses \ob{} and \cb, by `reading along the boundary' of $\la$, as follows. Working from left to right, for each node $\fka$ such that \nw\fka{} is not a node of $\la$, we write a \cb, and for each node $\fka$ such that \ne\fka{} is not a node of $\la$ we write a \ob.  (If neither \nw\fka{} nor \ne\fka{} is a node of $\la$, then we write \cb\ob.) We then append an infinite string of \ob s at the start of the sequence, and an infinite sequence of \cb s at the end.  The resulting doubly infinite sequence is called the \emph{parenthesis sequence} of $\la$.

\begin{eg}

Let $\la=(5,3^2,1)$. Then the parenthesis sequence of $\la$ is $\cdots\ob\ob\ob\cb\ob\cb\cb\ob\ob\cb\cb\ob\cb\cb\cdots$, as we see from the Young diagram, in which we mark the boundary in bold, extending it infinitely far to the north-west and north-east; segments $\diagdown$ contribute a \ob{} to the parenthesis sequence, while segments $\diagup$ contribute a \cb.
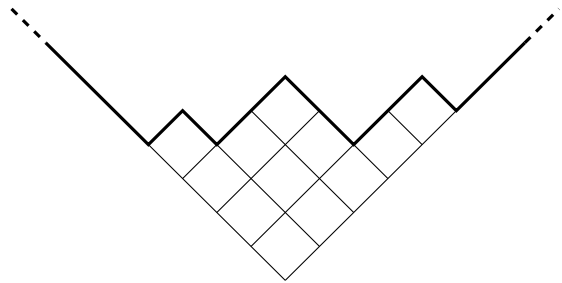
\[
\begin{tikzpicture}[scale=.45]
\foreach\x in{0,1}\draw(-\x,\x)--++(5,5);
\foreach\x in{2,3}\draw(-\x,\x)--++(3,3);
\draw(0,0)--(-4,4)--++(1,1)--++(4,-4);
\foreach\x in{2,3}\draw(\x,\x)--++(-3,3);
\foreach\x in{4,5}\draw(\x,\x)--++(-1,1);
\draw[very thick](-7,7)--++(3,-3)--++(1,1)--++(1,-1)--++(2,2)--++(2,-2)--++(2,2)--++(1,-1)--++(2,2);
\draw[very thick,dashed](-8,8)--(-7,7);
\draw[very thick,dashed](7,7)--(8,8);
\end{tikzpicture}
\]
\end{eg}

Now partition the parenthesis expression for $\la$ into pairs, in the usual way for pairing up parentheses: each \ob{} is paired with the first subsequent \cb{} for which there are an equal number of \cb s and \ob s in between.

\begin{eg}
Continuing the last example, we illustrate the pairs by numbering the $\cb$s in increasing order, and numbering each $\ob$ with the same number as its corresponding $\cb$:
\[
\begin{array}{cccccccccccccccccc}
\cdots&\ob&\ob&\ob&\ob&\cb&\ob&\cb&\cb&\ob&\ob&\cb&\cb&\ob&\cb&\cb&\cb&\cdots\\
\cdots&8&7&3&1&1&2&2&3&5&4&4&5&6&6&7&8&\cdots
\end{array}
\]
\end{eg}

\begin{rmk}
We remark that in \cite{kw} and elsewhere, finite parenthesis expressions are used: instead of appending infinite strings of \ob{} and \cb{} at the start and end of the expression obtained from the partition, one appends sufficiently long finite strings that the resulting expression is \emph{balanced}, meaning that it contains equally many \ob s and \cb s, and that any initial segment contains at least as many \ob s as \cb s. We find our convention more straightforward in the present context, and translation between the two conventions is very easy.
\end{rmk}

Now we recall a definition from \cite{kw}; it is phrased there in terms of finite Dyck paths, but it is more convenient for us to phrase it in terms of partitions. Given two partitions $\la,\mu$, write $\la\kwo\mu$ if the parenthesis expression for $\la$ can be obtained from that for $\mu$ by taking some of the pairs $\ob\cdots\cb$ and reversing them to get $\cb\cdots\ob$.

With this definition, we can describe the connection to \cedt s. The following \lcnamecref{cedtkw} is not hard to prove, and we leave it as an exercise for the reader.

\begin{propn}\label{cedtkw}
Suppose $\la,\mu$ are partitions. Then $\la\kwo\mu$ if and only if $\la\supseteq\mu$ and there is a \cedt{} of $\skw\la\mu$.
\end{propn}

Now we can describe the relationship between our results and those of Kenyon and Wilson.  Suppose $f,g$ are positive integers as in \cref{ypmsubsec}, and define matrices $N,P$ with rows and columns indexed by $\lfg$, and with
\[
N_{\la\mu}=(-1)^{|\la|+|\mu|}\et\la\mu,\qquad P_{\la\mu}=\ndt\la\mu.
\]
Then \cref{mainsn1,mainsn2} show that $N$ and $P$ are mutual inverses. In fact, Kenyon and Wilson prove this result, but in a slightly different form. They define a matrix $M$ with $M_{\la\mu}=1$ if $\la\kwo\mu$ and $0$ otherwise, and show (using a different recursion for $\ndt\la\mu$ from ours) that $(M^{-1})_{\la\mu}=(-1)^{|\la|+|\mu|}\ndt\la\mu$ \cite[Theorem 1.5]{kw}.  In view of \cref{cedtkw}, we have $M_{\la\mu}=\et\la\mu$, and so the main result of \cite{kw} shows that \cref{mainsn1,mainsn2} are equivalent. (There is another minor difference with our results, in that instead of the set $\lfg$ of partitions lying inside a rectangle, they consider the set of partitions lying inside a partition of the form $(c,c-1,\dots,1)$; however, one can easily show the equivalence of the two results by taking very large bounding partitions and appropriate submatrices.)

Of course, using this result we could remove a lot of the work in this paper, but we prefer to keep our proofs, thereby keeping our paper self-contained and providing a new (albeit longer) proof of the main result of \cite{kw}.

We remark that Kim \cite[\S7]{kim} observes a $q$-analogue of the above result (for which he also credits Konvalinka); in this, the integer $\ndt\la\mu$ is replaced with the polynomial
\[
\ndt\la\mu(q)=\sum_{T\in\dt\la\mu}q^{|T|},
\]
where $|T|$ denotes the number of tiles in the tiling $T$, and similarly for $\et\la\mu$. It should be possible to give a new proof of this result using our techniques: in our bijective results (\cref{noadd,bij1,bij1a,bij2,bij2a,bij3,exp1}) it is easy to keep track of how the number of tiles in a tiling changes, so $q$-analogues of our recursive results can easily be proved. Replacing the group algebra of the symmetric group with the Hecke algebra of type $A$ should enable $q$-analogues of \cref{mainsn1,mainsn2} to be proved, yielding a new proof of Kim's and Konvalinka's result. We leave the details to the proverbial interested reader.

\medskip

Another result from the literature relates to \cref{snsum}. To describe this relationship, we give an alternative characterisation of the function $\fp$, for which we need  another definition. Given a partition $\la$, construct the parenthesis sequence of $\la$ and divide it into pairs as above, and define a partial order on the set of pairs by saying that the pair $\underset a{\ob}\cdots\underset a{\cb}$ is larger than the pair $\underset b{\ob}\cdots\underset b{\cb}$ if the former is nested inside the latter: $\underset b\ob\cdots\underset{\vphantom ba}\ob\cdots\underset{\vphantom ba}\cb\cdots\underset b\cb$. Call this partially ordered set $\calp(\la)$.

\begin{eg}
Continuing the last example and retaining the labelling for the pairs in the parenthesis sequence for $\la$, we have the following Hasse diagram for $\calp(\la)$.
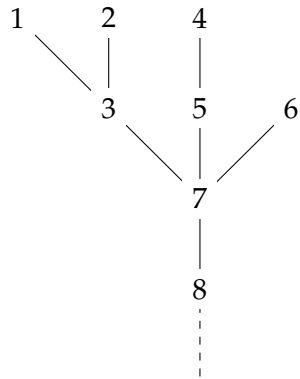
\[
\begin{tikzpicture}[scale=1.2]
\draw(0,0)--++(2,-2)--++(0,2);
\draw(1,0)--++(0,-1);
\draw(2,-2)--++(1,1);
\draw(2,-2)--++(0,-1);
\draw[dashed](2,-3)--++(0,-1);
\draw(0,0)node[fill=white]{1};
\draw(1,0)node[fill=white]{2};
\draw(2,0)node[fill=white]{4};
\draw(1,-1)node[fill=white]{3};
\draw(2,-1)node[fill=white]{5};
\draw(3,-1)node[fill=white]{6};
\draw(2,-2)node[fill=white]{7};
\draw(2,-3)node[fill=white]{8};
\end{tikzpicture}
\]
\end{eg}

Now we have the following.

\begin{propn}\label{fpchord}
Suppose $\la$ is a partition. Then $\fp(\la)$ is the number of linear extensions of the poset $\calp(\la)$.
\end{propn}

A proof of this \lcnamecref{fpchord} was sketched by David Speyer in a response to the author's MathOverflow question \cite{mo}. In a comment to the same answer, Philippe Nadeau pointed out that this provides a non-recursive expression for $\fp(\la)$. To state this, we let $\calp_N(\la)$ denote the poset obtained by taking the $N$ largest elements of $\calp(\la)$. Then for large enough $N$ the number of linear extensions of $\calp(\la)$ equals the number of linear extensions of $\calp_N(\la)$. Given a pair $p\in\calp(\la)$, define its \emph{length} $l(p)$ to be $1$ plus the number of intervening pairs. Now an easy exercise \cite[p.70, Exercise 20]{knuth} gives the following.

\begin{propn}\label{noext}
Suppose $\la$ is a partition, and $N\gg0$. Then the number of linear extensions of the poset $\calp_N(\la)$ equals
\[
\frac{N!}{\prod_{p\in\calp_N(\la)}l(p)}.
\]
\end{propn}

\begin{eg}
Taking $N=8$ in the last example and labelling each pair by its length, we have
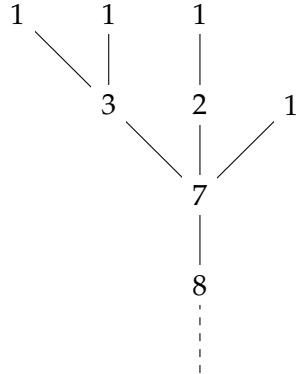
\[
\begin{tikzpicture}[scale=1.2]
\draw(0,0)--++(2,-2)--++(0,2);
\draw(1,0)--++(0,-1);
\draw(2,-2)--++(1,1);
\draw(2,-2)--++(0,-1);
\draw[dashed](2,-3)--++(0,-1);
\draw(0,0)node[fill=white]{1};
\draw(1,0)node[fill=white]{1};
\draw(2,0)node[fill=white]{1};
\draw(1,-1)node[fill=white]{3};
\draw(2,-1)node[fill=white]{2};
\draw(3,-1)node[fill=white]{1};
\draw(2,-2)node[fill=white]{7};
\draw(2,-3)node[fill=white]{8};
\end{tikzpicture}
\]
So $\fp(\la)$ is the number of linear extensions of $\calp_8(\la)$, which is
\[
\frac{8!}{2\times3\times7\times8}=120.
\]
\end{eg}

The author is very grateful to David Speyer and Philippe Nadeau for these comments (and also to Gjergji Zaimi for similar comments in the same thread), which inspired the introduction of Dyck tilings in the present work.

Given \cref{mainsn2} and \cref{fpchord}, \cref{snsum} is equivalent to the following theorem.

\begin{thm}\label{kimmain}
Suppose $\la$ is a partition and $N\gg0$. Then
\[
\sum_{\mu\subseteq\la}\ndt\la\mu=\frac{N!}{\prod_{p\in\calp_N(\la)}l(c)}.
\]
\end{thm}

This statement (in fact, a $q$-analogue) was conjectured by Kenyon and Wilson \cite[Conjecture 1]{kw}. It was proved inductively by Kim \cite[Theorem 1.1]{kim}, and then a bijective proof was given by Kim,  M\'esz\'aros, Panova and Wilson \cite[Theorem 1.1]{kmpw}. We retain our proof of \cref{snsum} to keep the paper self-contained (and so that \cref{snsum} can be proved without reference to tilings), and we note in passing that this yields a new proof of \cref{kimmain}.

\section{The homogeneous Garnir relations}\label{klrsec}

In this section we apply our earlier results to the study of representations of the (cyclotomic) quiver Hecke algebras (also known as \emph{KLR algebras}) of type $A$, and in particular the homogeneous Garnir relations for the \emph{universal graded Specht modules} of Kleshchev, Mathas and Ram \cite{kmr}.

\subsection{The quiver Hecke algebra}

Suppose $\bbf$ is a field of characteristic $p$, $e\in\{0,2,3,4,\dots\}$ and $n\in\bbn$. The \emph{quiver Hecke algebra} $\calr_n$ of type $A$ is a unital associative $\bbf$-algebra with a generating set
\[
\{y_1,\dots,y_n\}\cup\{\psi_1,\dots,\psi_{n-1}\}\cup\lset{e(\bfi)}{\bfi\in(\zez)^n}
\]
and a somewhat complicated set of defining relations, which may be found in \cite{bk} or \cite{kmr}, for example. These relations allow one to write down a basis for $\calr_n$: to do this, choose and fix a reduced expression $w=t_{i_1}\dots t_{i_r}$ for each $w\in\sss n$, and set $\psi_w=\psi_{i_1}\dots\psi_{i_r}$. Then the set
\[
\lset{\psi_wy_1^{c_1}\dots y_n^{c_n}e(\bfi)}{w\in\sss n,\ c_1,\dots,c_n\gs0,\ \bfi\in(\zez)^n}
\]
is a basis for $\calr_n$.

The quiver Hecke algebra has attracted considerable attention in recent years, thanks to the astonishing result of Brundan and Kleshchev \cite[Main Theorem]{bk} that when $p$ does not divide $e$, a certain finite-dimensional `cyclotomic' quotient of $\calr_n$ is isomorphic to the cyclotomic Hecke algebra of type $A$ (as introduced by Ariki--Koike \cite{arko} and Brou\'e--Malle \cite{brma}) defined at a primitive $e$th root of unity in $\bbf$; when $e=p$, there is a corresponding isomorphism to the degenerate cyclotomic Hecke algebra (which includes the group algebra $\bbf\sss n$ as a special case). This in particular shows that these Hecke algebras (which include the group algebra of the symmetric group) are non-trivially $\bbz$-graded, and has initiated the study of the graded representation theory of these algebras.

A crucial role in the representation theory of Hecke algebras is played by the \emph{Specht modules}. Brundan, Kleshchev and Wang \cite{bkw} showed how to work with the Specht modules in the quiver Hecke algebra setting, demonstrating in particular that these modules are graded. Kleshchev, Mathas and Ram \cite{kmr} gave a presentation for each Specht module with a single generator and a set of homogeneous relations. These relations include `homogeneous Garnir relations'; although these are in general simpler than the classical Garnir relations for ungraded Specht modules (which go back to \cite{garni} in the symmetric group case), the expressions given for these relations in \cite{kmr} are quite complicated. The purpose of this section is to use the results of the previous section to give a simpler expression for each homogeneous Garnir relation. In computations with graded Specht modules using the author's GAP programs, implementing these simpler expressions has been observed to have some benefits in terms of computational efficiency.

We now define the cyclotomic quotients of the quiver Hecke algebras and their Specht modules. Choose a positive integer $l$, and an $l$-tuple $(\kappa_1,\dots,\kappa_l)\in(\zez)^l$. For each $a\in\zez$ define $\kappa(a)$ to be the number of values of $j$ for which $\kappa_j=a$, and define $\calr^\kappa_n$ to be the quotient of $\calr_n$ by the ideal generated by the elements $y_1^{\kappa(i_1)}e(\bfi)$ for all $i\in(\zez)^n$. We use the same notation for the standard generators of $\calr_n$ and their images in $\calr^\kappa_n$.

\subsection{Row permutation modules and Specht modules}

From now on, we stick to the case $l=1$, which corresponds to the Iwahori--Hecke algebra of type $A$; there is no essential difference in the homogeneous Garnir relations for arbitrary $l$, and we save on notation and terminology by restricting to this special case. We also assume that $e>0$, since the difficulties we describe below do not arise in the case $e=0$. We work only with ungraded modules, since the grading plays no part in our results.

In the case $l=1$, the Specht modules (by which we mean the \emph{row Specht modules} of \cite[\S5.4]{kmr}) are labelled by partitions. For the sake of alignment with the rest of the literature, we now change to using the \emph{English} convention for Young diagrams in this section, where the first coordinate increases down the page and the second increases from left to right.

Suppose $\pi$ is a partition of $n$. If $(a,b)$ is a node, its \emph{residue} is defined to be $\res((a,b))=b-a+\kappa_1\ppmod e$. A \emph{$\pi$-tableau} is a bijection from the Young diagram of $\pi$ to the set $\{1,\dots,n\}$. We specify a particular $\pi$-tableau $\ttt^\pi$ by assigning the numbers $1,\dots,n$ to the nodes
\[
(1,1),(1,2),\dots,(1,\pi_1),\ (2,1),(2,2),\dots,(2,\pi_2),\ (3,1),(3,2),\dots
\]
in order. for example, when $\pi=(6,4,1^2)$ we have
\Yenglish
\[
\ttt^\pi=\young(123456,789\ten,\eleven,\twelve).
\]
We define an element $\bfi^\pi$ of $(\zez)^n$ by setting $\bfi^\pi_j=\res((\ttt^\pi)^{-1}(j))$. For example, with $\pi=(6,4,1^2)$, $e=4$ and $\kappa_1=0$, $\bfi^\pi$ is the sequence $(0,1,2,3,0,1,3,0,1,2,2,1)$.

Now we can define the \emph{row permutation module} $\yper\pi$ from \cite[\S5.3]{kmr}. This has a single generator $m^\pi$, and relations as follows:
\begin{itemize}
\item
$e(\bfi^\pi)m^\pi=m^\pi$;
\item
$y_jm^\pi=0$ for all $j$;
\item
$\psi_jm^\pi=0$ whenever $j$ and $j+1$ lie in the same row of $\ttt^\pi$.
\end{itemize}

It is easy to write down a basis for $\yper\pi$. Say that a $\pi$-tableau $\ttt$ is \emph{row-strict} if the entries increase from left to right along the rows. If $\ttt$ is row-strict, let $w^\ttt$ be the permutation taking $\ttt^\mu$ to $\ttt$, and define $\psi^\ttt:=\psi_{w^\ttt}$, and $m^\ttt=\psi^\ttt m^\pi$. Then by \cite[Theorem 5.6]{kmr} the set $\lset{m^\ttt}{\ttt\text{ a row-strict $\pi$-tableau}}$ is a basis for $\yper\pi$.

The Specht module $\spe\pi$ is the quotient of $\yper\pi$ by the homogeneous Garnir relations, which we now define. Say that a node $\fkn=(a,b)$ of $\pi$ is a \emph{Garnir node} if $(a+1,b)$ is also a node of $\pi$. If $\fkn$ is a Garnir node, define the \emph{Garnir belt} $\blt^{\fkn}$ to be the set of nodes
\[
\{(a,b),(a,b+1),\dots,(a,\pi_a),\ (a+1,1),(a+1,2),\dots,(a+1,b)\}.
\]
The \emph{Garnir tableau} $\ttg^{\fkn}$ is the tableau defined by taking $\ttt^\pi$ and rearranging the entries within $\blt^{\fkn}$ so that they increase from bottom left to top right. For example, with $\pi=(6,4,1^2)$ and $\fkn=(1,3)$, the tableau $\ttg^{\fkn}$ (with the Garnir belt $\blt^\fkn$ shaded) is as follows.
\[
\begin{tikzpicture}[scale=1]
\Yfillcolour{white!80!black}
\tgyoung(0cm,0cm,::;;;;,;;;)
\Yfillcolour{white}
\Yfillopacity0
\tyoung(0cm,0cm,126789,345\ten,\eleven,\twelve)
\end{tikzpicture}
\]
Now we define \emph{bricks}: an $\fkn$-brick is defined to be a set of $e$ consecutive nodes in the same row of $\blt^{\fkn}$ of which the leftmost has the same residue as $\fkn$. For example, take $\pi=(6,4,1^2)$ and $\fkn=(1,3)$ as above. If $e=2$, then the bricks are
\[
\{(1,3),(1,4)\},\{(1,5),(1,6)\},\{(2,2),(2,3)\},
\]
while if $e=3$ the bricks are
\[
\{(1,3),(1,4),(1,5)\},\{(2,1),(2,2),(2,3)\}.
\]
Let $f$ denote the number of bricks in row $a$ and $g$ the number of bricks in row $a+1$, and set $k=f+g$. If $k>0$, let $d=d^{\fkn}$ be the smallest number which is contained in a brick in $\ttg^\fkn$, and for $1\ls i<k$ define $w_i$ to be the permutation
\[
(d+ie-e,\ d+ie)(d+ie-e+1,\ d+ie+1),\dots,(d+ie-1,\ d+ie+e-1),
\]
Recalling the notation $\psi_w$ for $w\in\sss n$ from above, set $\sigma^\fkn_i=\psi_{w_i}$; unlike the situation for arbitrary $w$, $\psi_{w_i}$ does not depend on the choice of reduced expression for $w_i$. Also set $\tau^\fkn_i=\sigma^\fkn_i+1$.

\begin{rmk}
In the definition \cite[(5.7)]{kmr}, the elements $\tau^\fkn_i,\sigma^\fkn_i$ include an idempotent $e(\bfi^\fkn)$ as a factor; but this factor is unnecessary for our purposes, so we prefer the simpler version.
\end{rmk}

Now let $\ttt^{\fkn}$ be the tableau obtained from $\ttg^{\fkn}$ by re-ordering the bricks so that their entries increase along row $a$ and then row $a+1$.

Continuing the example above, we illustrate the tableau $\ttt^\fkn$ for $e=2,3$, shading the Garnir belt and outlining the bricks:
\[
\ttt^{\fkn}=
\begin{tikzpicture}[scale=1,baseline=-17pt]
\Yfillcolour{white!80!black}
\Yfillopacity1
\tgyoung(0cm,0cm,::;;;;,;;;)
\Yfillopacity0
\tyoung(0cm,0cm,124567,389\ten,\eleven,\twelve)
\Ylinethick{1.3pt}
\tgyoung(0cm,0cm,::_2_2,:_2\ )
\Ylinethick{.3pt}
\end{tikzpicture}
\text{ if }e=2,\qquad
\ttt^{\fkn}=
\begin{tikzpicture}[scale=1,baseline=-17pt]
\Yfillcolour{white!80!black}
\Yfillopacity1
\tgyoung(0cm,0cm,::;;;;,;;;)
\Yfillopacity0
\tyoung(0cm,0cm,123459,678\ten,\eleven,\twelve)
\Ylinethick{1.3pt}
\tgyoung(0cm,0cm,::_3,_3\ )
\Ylinethick{.3pt}
\end{tikzpicture}
\text{ if }e=3.
\]
Define $\psi^{\ttt^\fkn}$ as above; as noted in \cite[\S5.4]{kmr}, $\psi^{\ttt^\fkn}$ is independent of the choice of reduced expression for $w_{\ttt^\fkn}$.

Now we define elements of $\calr^\kappa_n$ corresponding to the elements $t_\la$ and $s_\la$ of $\bbf\sss k$ defined in \cref{ypmsec}. For $\la\in\lfg$, define
\[
\tau^\fkn_\la=\left(\tau^\fkn_{\la_1+g}\tau^\fkn_{\la_1+g+1}\dots \tau^\fkn_{f+g-1}\right)\left(\tau^\fkn_{\la_2+g-1}\tau^\fkn_{\la_2+g}\dots \tau^\fkn_{f+g-2}\right)\dots\left(\tau^\fkn_{\la_g+1}\tau^\fkn_{\la_g+2}\dots \tau^\fkn_f\right),
\]
and define $\sigma^\fkn_\la$ similarly. Then the elements $\sigma^\fkn_\la\psi^{\ttt^\fkn}$ are precisely the elements $\psi^\ttt$ for row-strict tableaux $\ttt$ obtained from $\ttg^\fkn$ by permuting the bricks.

Now we can define the homogeneous Garnir relations and the Specht module. Given a Garnir node $\fkn$, the corresponding \emph{Garnir element} is
\[
g^\fkn=\zsum{\la\in\lfg}\tau^\fkn_\la\psi^{\ttt^\fkn}.
\]
The Specht module $\spe\pi$ is defined to be the quotient of $\yper\pi$ by the relations $g^\fkn m^\pi=0$ for all Garnir nodes $\fkn$ of $\pi$.

\subsection{Re-writing the Garnir relations}

The difficulty with the definition of the Specht module given above, especially from a computational point of view, is that the individual terms $\tau^\fkn_\la\psi^{\ttt^\fkn}m^\pi$ appearing in the Garnir relation are not readily expressed in terms of the standard basis $\lset{m^\ttt}{\ttt\text{ row-strict}}$ for $\yper\pi$; to express it in these terms, some quite involved reduction is required using the defining relations in $\calr_n$.

\begin{eg}
Take $\pi=(11,5,3,1)$, $\fkn=(1,5)$ and $e=3$. Then we have $f=2$, $g=1$, and
\[
\sigma^\fkn_1=\psi_9\psi_8\psi_7\psi_{10}\psi_9\psi_8\psi_{11}\psi_{10}\psi_9,\qquad\sigma^\fkn_2=\psi_{12}\psi_{11}\psi_{10}\psi_{13}\psi_{12}\psi_{11}\psi_{14}\psi_{13}\psi_{12},
\]
while $\psi^{\ttt^\fkn}=\psi_{6}\psi_{5}\psi_{7}\psi_{6}\psi_{8}\psi_{7}\psi_{9}\psi_{8}\psi_{10}\psi_{9}\psi_{11}\psi_{10}\psi_{15}\psi_{14}\psi_{13}\psi_{12}\psi_{11}$. Hence the Garnir element $g^\fkn$ is
\begin{align*}
\left(\tau^\fkn_\varnothing+\tau^\fkn_{(1)}+\tau^\fkn_{(2)}\right)\psi^{\ttt^\fkn}&=\left(\left(\sigma^\fkn_1+1\right)\left(\sigma^\fkn_2+1\right)+\left(\sigma^\fkn_2+1\right)+1\right)\psi^{\ttt^\fkn}\\
&=\left(\sigma^\fkn_1\sigma^\fkn_2+\sigma^\fkn_1+2\sigma^\fkn_2+3\right)\psi^{\ttt^\fkn}.
\end{align*}
A non-trivial calculation using the defining relations for $\calr_n$ and the fact that $m^\pi$ is killed by $\psi_5$, $\psi_6$ and $\psi_7$ shows that $\sigma^\fkn_1\psi^{\ttt^\fkn}m^\pi=0$. Hence
\[
g^\fkn m^\pi=\left(\sigma^\fkn_\varnothing+2\sigma^\fkn_{(1)}+3\sigma^\fkn_{(2)}\right)\psi^{\ttt^\fkn}m^\pi=m^{\ttg^\fkn}+2m^\ttu+3m^{\ttt^\fkn},
\]
where
\begin{align*}
\ttg^\fkn&=\young(1234\ten\eleven\twelve\thirteen\fourteen\fifteen\sixteen,56789,\seventeen\eighteen\nineteen,\twenty),\\
\ttu&=\young(1234789\thirteen\fourteen\fifteen\sixteen,56\ten\eleven\twelve,\seventeen\eighteen\nineteen,\twenty),\\
\ttt^\fkn&=\young(1234789\ten\eleven\twelve\sixteen,56\thirteen\fourteen\fifteen,\seventeen\eighteen\nineteen,\twenty).
\end{align*}
\end{eg}

We wish to generalise the above example, to re-write $g^\fkn m^\pi$ for an arbitrary Garnir node $\fkn$ as a linear combination of elements $\sigma^\fkn_\la\psi^{\ttt^\fkn}m^\pi$. Fortunately, all the computations that we need with the defining relations for $\calr_n$ have already been done in \cite{kmr}.

Given a Garnir node $\fkn$ of $\pi$, define the \emph{brick permutation space} $T^{\pi,\fkn}$ to be the $\bbf$-subspace of $\yper\pi$ spanned by all elements of the form $\sigma^\fkn_{i_1}\dots\sigma^\fkn_{i_s}\psi^{\ttt^\fkn}m^\pi$. Then the $\sigma^\fkn_i$, and hence the $\tau^\fkn_i$, act on $T^{\pi,\fkn}$, and we have the following.

\begin{thmc}{kmr}{Theorem 5.11}\label{brickspace}
As operators on $T^{\pi,\fkn}$, the elements $\tau^\fkn_1,\dots,\tau^\fkn_{k-1}$ satisfy the Coxeter relations for the symmetric group $\sss k$, and hence $T^{\pi,\fkn}$ can be considered as an $\bbf\sss k$-module. In fact, as an $\bbf\sss n$-module $T^{\pi,\fkn}$ is isomorphic to the Young permutation module $\mfg$, with an isomorphism given by mapping $\psi^{\ttt^\fkn}m^\pi$ to the standard generator $\gn$.
\end{thmc}

By \cref{brickspace} and the discussion at the start of \cref{ypmsec}, the sets
\[
\lset{\tau^\fkn_\la\psi^{\ttt^\fkn}m^\pi}{\la\in\lfg}\qquad\text{and}\qquad\lset{\sigma^\fkn_\la\psi^{\ttt^\fkn}m^\pi}{\la\in\lfg}
\]
are $\bbf$-bases of $T^{\pi,\fkn}$. \cref{mainsn1,mainsn2} give the change-of-basis matrices between these bases, while \cref{snsum} gives the sum of the elements of the first basis in terms of the second basis. So we can immediately deduce the following, which is the main result of this section.

\begin{thm}\label{maingarn}
Suppose $\pi$ is a partition and $\fkn=(a,b)$ is a Garnir node of $\pi$, and let $f,g$ denote the numbers of bricks in rows $a$ and $a+1$ of $\blt^\fkn$ respectively. Define $\sigma^\fkn_\la$ and $\tau^\fkn_\la$ for $\la\in\lfg$ as above, and define the \emph{modified Garnir element}
\[
\hat g^\fkn=\zsum{\la\in\lfg}\fp(\la)\sigma^\fkn_\la\psi^{\ttt^\fkn}.
\]
Then:
\begin{enumerate}
\item
$\sigma^\fkn_\mu\psi^{\ttt^\fkn}m^\pi=\zsum{\la\in\lfg}(-1)^{|\la|+|\mu|}\et\la\mu\tau^\fkn_\la\psi^{\ttt^\fkn}m^\pi$;
\item
$\tau^\fkn_\mu\psi^{\ttt^\fkn}m^\pi=\zsum{\la\in\lfg}\ndt\la\mu\sigma^\fkn_\la\psi^{\ttt^\fkn}m^\pi$;
\item
$\hat g^\fkn m^\pi=g^\fkn m^\pi$.
\end{enumerate}
Hence the Specht module $\spe\pi$ is the quotient of $\yper\pi$ by the relations $\hat g^\fkn m^\pi=0$ for all Garnir nodes $\fkn$ of $\pi$.
\end{thm}

\begin{rmks}\indent
\begin{enumerate}
\vspace{-\topsep}
\item
Statements (1) and (2) in the theorem are not necessary for the Garnir relations, but may be of interest for computation in the row permutation modules; indeed, the coefficients in (1) are essentially the coefficients $c_w$ appearing in \cite[Corollary 5.12]{kmr}.
\item
We emphasise that the identities in \cref{maingarn} are not true if the terms $m^\pi$ are omitted; in general, $\tau^\fkn_\mu$ is not a linear combination of $\sigma^\fkn_\la$s (even if the idempotent $e(\bfi^\fkn)$ from \cite[(5.7)]{kmr} is reinstated in the definition of $\sigma^\fkn_i$ and $\tau^\fkn_i$).
\end{enumerate}
\end{rmks}

We end with an example which illustrates how the Garnir relations are used in Specht module computations. Say that a $\pi$-tableau is \emph{standard} if the entries are increasing both along the rows and down the columns. If we let $v^\ttt$ denote the image of $m^\ttt$ in the Specht module, then by \cite[Corollary 6.24]{kmr} the set $\lset{v^\ttt}{\ttt\text{ a standard $\pi$-tableau}}$ is a basis for $\spe\pi$. The Garnir relations allow one to write $v^\ttt$ as a linear combination of this basis when $\ttt$ is a row-strict, but not standard, tableau. The foundation for this is \cite[Lemma 5.13]{kmr} which does this in the case where $\ttt=\ttg^\fkn$ for a Garnir node $\fkn$.

For our example, take $\pi=(8,4)$ and $e=2$. Let $\fkn=(1,4)$. Then $f=g=2$, and
\[
\ttg^\fkn=\young(12389\ten\eleven\twelve,4567).
\]
As with any Garnir node, $\ttg^\fkn$ is row-strict but not standard. The modified Garnir element is
\[
\hat g^\fkn=\left(\sigma^\fkn_\varnothing+2\sigma^\fkn_{(1)}+3\sigma^\fkn_{(2)}+3\sigma^\fkn_{(1^2)}+6\sigma^\fkn_{(2,1)}+6\sigma^\fkn_{(2^2)}\right)\psi^{\ttt^\fkn},
\]
and the corresponding Garnir relation yields
\[
v^{\ttg^\fkn}=-2v^\ttu-3v^\ttv-3v^\ttw-6v^\ttx-6v^\tty,
\]
where $\ttu,\ttv,\ttw,\ttx,\tty$ are standard tableaux that we describe as follows. We represent the elements $\psi^\ttu,\psi^\ttv,\psi^\ttw,\psi^\ttx,\psi^\tty$ using the braid diagrams of Khovanov and Lauda \cite{kl}; by using shaded squares, we illustrate (reverting to the Russian convention!) how the labelling partitions $\la\in\scrp_{2,2}$ arise.
{\allowdisplaybreaks
\begin{alignat*}3
\ttg^\fkn&=\young(12389\ten\eleven\twelve,4567),&\qquad\psi^{\ttg^\fkn}&=\sigma^\fkn_\varnothing\psi^{\ttt^\fkn}&&=\ 
\begin{tikzpicture}[scale=.3,baseline=1.7cm,rounded corners=2pt]
\foreach\x in{1,2,3}\draw(\x,0)--++(0,9);
\foreach\x in{0,1,2,3}\draw(4+\x,0)--++(0,3.5-\x)--++(5,5)--++(0,.5+\x);
\foreach\x in{0,1,2,3,4}\draw(8+\x,0)--++(0,.5+\x)--++(-4,4)--++(0,4.5-\x);
\end{tikzpicture}
\\
\ttu&=\young(12367\ten\eleven\twelve,4589),&\psi^{\ttu}&=\sigma^\fkn_{(1)}\psi^{\ttt^\fkn}&&=\ 
\begin{tikzpicture}[scale=.3,baseline=1.7cm]
\draw[fill=white!85!black](7.5,0)--++(2,2)--++(-2,2)--++(-2,-2);
\foreach\x in{0,2}{\draw[line width=2pt,white](7.5+\x,\x)--++(-2,2);\draw[line width=2pt,white](7.5-\x,\x)--++(2,2);}
\foreach\x in{1,2,3}\draw[rounded corners](\x,0)--++(0,9);
\foreach\x in{0,1}\draw[rounded corners](4+\x,0)--++(0,3.5-\x)--++(5,5)--++(0,.5+\x);
\foreach\x in{2,3,4}\draw[rounded corners](8+\x,0)--++(0,.5+\x)--++(-4,4)--++(0,4.5-\x);
\foreach\x in{0,1}\draw[rounded corners](6+\x,0)--++(0,2.5+\x)--++(-2,2)--++(0,4.5-\x);
\foreach\x in{0,1}\draw[rounded corners](8+\x,0)--++(0,3.5-\x)--++(3,3)--++(0,2.5+\x);
\end{tikzpicture}
\\
\ttv&=\young(1236789\twelve,45\ten\eleven),&\psi^{\ttv}&=\sigma^\fkn_{(2)}\psi^{\ttt^\fkn}&&=\ 
\begin{tikzpicture}[scale=.3,baseline=1.7cm]
\draw[fill=white!85!black](7.5,0)--++(4,4)--++(-2,2)--++(-4,-4);
\foreach\x in{0,2,4}\draw[line width=2pt,white](7.5+\x,\x)--++(-2,2);
\foreach\x in{0,2}\draw[line width=2pt,white](7.5-\x,\x)--++(4,4);
\foreach\x in{1,2,3}\draw[rounded corners](\x,0)--++(0,9);
\foreach\x in{0,1}\draw[rounded corners](4+\x,0)--++(0,3.5-\x)--++(5,5)--++(0,.5+\x);
\foreach\x in{4}\draw[rounded corners](8+\x,0)--++(0,.5+\x)--++(-4,4)--++(0,4.5-\x);
\foreach\x in{0,1,2,3}\draw[rounded corners](6+\x,0)--++(0,2.5+\x)--++(-2,2)--++(0,4.5-\x);
\foreach\x in{0,1}\draw[rounded corners](10+\x,0)--++(0,5.5-\x)--++(1,1)--++(0,2.5+\x);
\end{tikzpicture}
\\
\ttw&=\young(12345\ten\eleven\twelve,6789),&\psi^{\ttw}&=\sigma^\fkn_{(1^2)}\psi^{\ttt^\fkn}&&=\ 
\begin{tikzpicture}[scale=.3,baseline=1.7cm]
\draw[fill=white!85!black](7.5,0)--++(2,2)--++(-4,4)--++(-2,-2);
\foreach\x in{0,2}\draw[line width=2pt,white](7.5+\x,\x)--++(-4,4);
\foreach\x in{0,2,4}\draw[line width=2pt,white](7.5-\x,\x)--++(2,2);
\foreach\x in{1,2,3,4,5}\draw[rounded corners](\x,0)--++(0,9);
\foreach\x in{0,1}\draw[rounded corners](6+\x,0)--++(0,5.5-\x)--++(3,3)--++(0,.5+\x);
\foreach\x in{2,3,4}\draw[rounded corners](8+\x,0)--++(0,.5+\x)--++(-4,4)--++(0,4.5-\x);
\foreach\x in{0,1}\draw[rounded corners](8+\x,0)--++(0,3.5-\x)--++(3,3)--++(0,2.5+\x);
\end{tikzpicture}
\\
\ttx&=\young(1234589\twelve,67\ten\eleven),&\psi^{\ttx}&=\sigma^\fkn_{(2,1)}\psi^{\ttt^\fkn}&&=\ 
\begin{tikzpicture}[scale=.3,baseline=1.7cm]
\draw[fill=white!85!black](7.5,0)--++(4,4)--++(-2,2)--++(-2,-2)--++(-2,2)--++(-2,-2);
\foreach\x in{0,2}{\draw[line width=2pt,white](7.5+\x,\x)--++(-4,4);\draw[line width=2pt,white](7.5-\x,\x)--++(4,4);}
\foreach\x in{4}{\draw[line width=2pt,white](7.5+\x,\x)--++(-2,2);\draw[line width=2pt,white](7.5-\x,\x)--++(2,2);}
\foreach\x in{1,2,3,4,5}\draw[rounded corners](\x,0)--++(0,9);
\foreach\x in{0,1}\draw[rounded corners](6+\x,0)--++(0,5.5-\x)--++(3,3)--++(0,.5+\x);
\foreach\x in{4}\draw[rounded corners](8+\x,0)--++(0,.5+\x)--++(-4,4)--++(0,4.5-\x);
\foreach\x in{0,1}\draw[rounded corners](10+\x,0)--++(0,5.5-\x)--++(1,1)--++(0,2.5+\x);
\foreach\x in{0,1}\draw[rounded corners](8+\x,0)--++(0,4.5+\x)--++(-2,2)--++(0,2.5-\x);
\end{tikzpicture}
\\
\tty&=\young(1234567\twelve,89\ten\eleven),&\psi^{\tty}&=\sigma^\fkn_{(2^2)}\psi^{\ttt^\fkn}&&=\ 
\begin{tikzpicture}[scale=.3,baseline=1.7cm]
\draw[fill=white!85!black](7.5,0)--++(4,4)--++(-4,4)--++(-4,-4);
\foreach\x in{0,2,4}{\draw[line width=2pt,white](7.5+\x,\x)--++(-4,4);\draw[line width=2pt,white](7.5-\x,\x)--++(4,4);}
\foreach\x in{1,2,3,4,5,6,7}\draw[rounded corners](\x,0)--++(0,9);
\foreach\x in{4}\draw[rounded corners](8+\x,0)--++(0,.5+\x)--++(-4,4)--++(0,4.5-\x);
\foreach\x in{-2,-1,0,1}\draw[rounded corners](10+\x,0)--++(0,5.5-\x)--++(1,1)--++(0,2.5+\x);
\end{tikzpicture}
\end{alignat*}}

\end{document}